\documentclass[11pt]{amsart}
\input{xy}
\usepackage{amscd}
\xyoption{all}

\usepackage{amsfonts}
\usepackage{amssymb}

\title[Quivers and adjoint action]{Adjoint action of automorphism groups on
radical endomorphisms, generic equivalence and Dynkin quivers
}

\author{Bernt Tore Jensen and Xiuping Su}

\newtheorem{theorem}{{Theorem}}

\newtheorem{lemma}{{Lemma}}
\newtheorem{definition}[]{{Definition}}
\newtheorem{remark}{{Remark}}[section]
\newtheorem{corollary}[lemma]{{Corollary}}
\newtheorem{example}[]{{Example}}

\newcommand{\lra}{\longrightarrow}

\newcommand{\N}{{\mathbb N}}

\newcommand{\Aut}{\mathrm{Aut}}
\newcommand{\radE}{\mathrm{radEnd}}
\newcommand{\radH}{\mathrm{radHom}}
\newcommand{\rad}{\mathrm{rad}}
\newcommand{\Ext}{\mathrm{Ext}}
\newcommand{\Rep}{\mathrm{Rep}}
\newcommand{\Hom}{\mathrm{Hom}}
\newcommand{\dimv}{\underline{\mathrm{dim}}}
\newcommand{\End}{\mathrm{End}}
\newcommand{\Gl}{\mathrm{Gl}}
\newcommand{\Supp}{\mathrm{Supp}}
\renewcommand{\Im}{\mathrm{Im}}
\newcommand{\Ker}{\mathrm{Ker}}
\newcommand{\Mat}{\mathrm{Mat}}
\newcommand{\B}{\mathrm{B}}

 \normalbaselines

\newcommand{\ra}{\rightarrow}
\newcommand{\sdp}{\times\kern-.2em\vrule height1.1ex depth-.05ex}
\newcommand{\epi}{\lra \kern-.8em\ra}

\begin{document}

\maketitle

\begin{abstract}
Let $Q$ be a connected quiver with no oriented cycles, $k$ the field of
complex numbers and $P$ a projective 
representation of $Q$. We study the adjoint action of the automorphism 
group $\Aut_{kQ} P$ on the space of radical endomorphisms $\radE_{kQ}P$.  
Using generic equivalence, we show that the quiver $Q$ has the 
property that there exists a dense open 
$\Aut_{kQ} P$-orbit in $\radE_{kQ} P$, for all projective representations $P$, if 
and only if $Q$ is a Dynkin quiver. This gives a new characterisation of 
Dynkin quivers.
\end{abstract}

\section*{Introduction}

Let $A=kQ$ be the path algebra of a connected quiver $Q$ with no oriented cycles,
and $k$ is the field of complex numbers. Let $P$ be a projective representation
of $Q$. Note that there is an adjoint action of $\Aut_A P$ on the Jacobson radical
$\radE_A P$ of the endomorphism ring $\End_A P$, i.e. $g\cdot f = gfg^{-1}$ for 
any $g\in \Aut_A P$ and $f\in \radE_A P$. We are interested in generic 
$\Aut_A P$-orbits in $\radE_A P$, in particular in the existence of open orbits.

If $Q$ is of type $\mathbb{A}$ with linear orientation, then $\End_A P$ is a 
parabolic Lie algebra and $\radE_A P$ has a dense open $\Aut_A P$-orbit by a 
theorem of Richardson \cite{richardson}. The theorem of Richardson holds for
parabolic Lie algebras of any reductive Lie-algebra, and in type $\mathbb{A}$,
elements with dense orbits have been explicitly constructed by Br\"{u}stle, Hille, 
Ringel and R\"{o}rhle using representations of quivers \cite{BHRR}, and in the
classical types by Baur using a different approach \cite{baur}. These results and
methods have been extended in various directions, see for example 
\cite{erdmann, hille1, hille2, JSY, tan}.

If $Q$ is of type $\mathbb{A}$ with arbitrary orientation, then $\End_A P$ is a 
seaweed Lie algebra \cite{S1,S2} (also called a biparabolic Lie algebra \cite{B1}), 
and $\radE_A P$ has 
a dense open $\Aut_A P$-orbit by a theorem of Jensen, Su and Yu \cite{JSY}. At present 
it is not known if seaweed Lie algebras of other classical types have dense orbits,
however a counterexample exists for Lie algebras of type $\mathbb{E}_8$ 
\cite{JSY}. 

If $Q$ is of infinite representation type, i.e. $Q$ is not Dynkin, then there is not necessarily 
a dense open $\Aut_A P$-orbit in $\radE_A P$.

\begin{lemma} \label{nonfinite}
Let $Q$ be a quiver which is not of Dynkin type. Then there exists a projective 
representation $P$ of $Q$, such that there is no dense open $\Aut_A P$-orbit in 
$\radE_A P$.
\end{lemma}

The main result of this paper is as follows.

\begin{theorem} \label{QuiverTheorem}
Let $Q$ be a Dynkin-quiver and let $P$ be an arbitrary projective representation of 
$Q$. Then there is a unique dense open $\Aut_A P$-orbit in $\radE_A P$. 
\end{theorem}

By Gabriel's Theorem \cite{Gabriel}, Dynkin quivers are exactly those quivers 
with finite representation type. Lemma \ref{nonfinite} and Theorem 
\ref{QuiverTheorem} together give a new characterisation of Dynkin quivers, that is, 
a quiver is Dynkin if and only if there is a dense open $\Aut_A P$-orbit in $\radE_A P$ 
for all projective representations $P$ of $Q$.  

We emphasize that although there is a dense open $\Aut_A P$-orbit in $\radE_A P$, 
when $Q$ is Dynkin, there are in general infinitely many $\Aut_A P$-orbits in $\radE_A P$. 
This is an interesting phenomenon that can already be seen in type $\mathbb{A}$ 
(see \cite{DR, BHRR}). 

Note that the main result in \cite{HV} connects the adjoint action of $\Aut_AP$ on 
$\radE_A P$ with good representations of a double quiver of $Q$ with relations. 
In particular, there is an open dense $\Aut_AP$-orbit in $\radE_AP$ if and only if 
there is an open orbit in the variety of good representations. To prove our main result, 
we show that for some quivers $Q$, including Dynkin quivers other than type $\mathbb{A}$, 
the varieties of good representations are generically equivalent to representation 
varieties of a quiver $Q'$ with the same underlying graph as $Q$. This means that, although 
the geometric properties of the varieties of good representations are in general very
complicated, the behavior of generic good representations is simpler, and 
similar to that of generic representations of $Q'$. In particular, there is a dense 
open orbit in the variety of good representations if and only if there is an open 
orbit in the corresponding representation variety of $Q'$.

The paper is organized as follows. In Section 1 we recall basic 
facts on quivers and their representations, and introduce the notion of a
generic section and generic equivalence. In Section 2 we recall results of 
Br\"{u}stle and Hille \cite{HB}, and Hille and Vossieck \cite{HV} on the use of double 
quivers to parameterise $\Aut_A P$-orbits of radical endomorphisms of projective 
representations $P$. In Section 3 we recall Voigt's lemma for quiver representations 
and give a criterion for rigidity of good representations. 
In Section 4, using the criterion for rigidity, we give new simpler proofs to show that 
representations constructed for type 
$\mathbb{A}$ in \cite{BHRR, JSY} have dense open orbits. Using this construction 
in type $\mathbb{A}$ and the criterion in Section 3, we prove the main results in 
Section 5 and Section 6.

\section{Representation varieties and generic sections}

\subsection{Representation varieties}

Let $Q$ be a quiver with $Q_0=\{1,\cdots,n\}$ the set of vertices and $Q_1$ the 
set of arrows. Let $s,t:Q_1\rightarrow Q_0$ be the functions mapping an arrow to 
its starting and terminating vertex, respectively. A vertex $i\in Q_0$ is called a sink if 
there are no arrows starting at $i$, and a source if there are no arrows terminating 
at $i$. It is called admissible if it is either a sink or a source, and interior if there 
are at least two arrows incident to $i$.

A representation $M$ of $Q$ consists of vector spaces $\{M_i\}_{i\in Q_0}$ and 
linear maps $\{M_\alpha:M_{s(\alpha)}\rightarrow M_{t(\alpha)}\}_{
\alpha\in Q_1}$. A homomorphism of representations $h:M\rightarrow N$ is a 
collection of maps $h_i:M_i\rightarrow N_i$ satisfying 
$h_jM_\alpha=N_\alpha h_i$ for each arrow $\alpha:i\rightarrow j\in Q_1$. 
The direct sum of two representations is obtained by taking direct sums of vector 
spaces at each vertex and direct sum of linear maps at each arrow. A 
representation is indecomposable if it is not isomorphic to the direct sum of two 
nonzero representations.

For a representation $M$, let $\dimv M=(\dim_k M_i)_{i\in Q_0}$ denote the 
dimension vector of $M$. Let $c\in \mathbb{N}^n$ be a dimension 
vector and 
$$
\Rep (Q,c)=\prod_{\alpha\in Q_1}\Hom_k(k^{c_{s(\alpha)}}, k^{c_{t(\alpha)}})
$$ 
be the space of representations. We fix a basis and view elements in $\Rep (Q,c)$ as tuples 
of matrices. Let $\Gl_{c_i}$ denote the general
linear group of invertible $c_i\times c_i$-matrices. The group 
$$\Gl(c)=\prod_{i\in Q_0}\Gl_{c_i}$$ 
acts on $\Rep (Q,c)$ by change of basis. There is a bijection between $\Gl(c)$-orbits 
in $\Rep (Q,c)$ and isomorphism classes of representations of $Q$ with dimension 
vector $c$. 

The path algebra $A=kQ$ is the algebra with basis the set of paths in $Q$. For 
two paths $p$ and $q$, their product is defined to be the composition $pq$, if $q$ 
terminates where $p$ starts, and zero otherwise. For each vertex $i\in Q_0$, let 
$e_i$ denote the trivial path of length zero at $i$. The trivial paths form a set of 
pairwise orthogonal idempotents for $A$. There is an equivalence of categories 
between left $A$-modules and representations of $Q$. Using this equivalence,
any representation of $Q$ can be viewed as an $A$-module, and vice versa.

A quiver $Q$ is Dynkin if the underlying graph of $Q$ is one of the Dynkin graphs 
$\mathbb{A}_i,\mathbb{D}_j,\mathbb{E}_l$ for $i\geq 1, j\geq 4, l=6,7,8$. 
If $Q$ is Dynkin, by Gabriel's Theorem \cite{Gabriel} there are only finitely 
many orbits in $\Rep (Q,c)$ for any dimension vector $c$ and thus there is always 
a dense orbit in $\Rep (Q,c)$. Moreover, the Dynkin quivers are characterised by 
this property. We summarise these properties in the following theorem.

\begin{theorem}[\cite{Gabriel}] \label{Gabriel}
If $Q$ is a Dynkin quiver and $c$ is a dimension vector, then there is a dense 
open $\Gl(c)$-orbit in $\Rep (Q,c)$. Moreover, if $Q$ is not a Dynkin quiver, 
then there is a dimension vector $c$ such that there is no dense open 
$\Gl(c)$-orbit in $\Rep (Q,c)$.
\end{theorem}

We also need representations satisfying relations. Let $I\subseteq A$ be an ideal.
The corresponding subset $$\Rep (A/I,c)\subseteq \Rep (Q,c),$$ 
consisting of representations that are annihilated by $I$ is a $\Gl(c)$-stable Zariski closed 
subvariety, which is called the {\it variety of $(Q,I)$-representations} 
with dimension vector $c$.

\subsection{Generic sections}\label{relsec}

We will relate $\Aut_A P$-orbits in $\radE_A P$ to representations 
of a quiver without relations. This relationship will be made precise
using generic equivalence to be defined below.

A $G$-space is an affine space $V$ with a regular $G$-action for a 
connected algebraic group $G$.


\begin{definition} \label{gensectiondef}
An $H$-space $W$ is called a {\it generic section} of the $G$-space $V$
if there is an injective morphism $\phi:W\lra V$ such that 
\begin{itemize}
\item[(1)] $\Im\phi\subseteq V$ is an affine subspace. 
\item[(2)] $G\cdot \Im\phi$ contains a non-empty open subset of $V$.
\item[(3)] There is a non-empty open subset $W'\subseteq W$ such that
$\phi(H\cdot w)=(G\cdot \phi(w))\cap \Im\phi$ for all $w\in W'$.
\end{itemize}
\end{definition}

We give two examples of generic sections, where the first one will be
used in the proof of Lemma \ref{l1} in Section 6.

\begin{example} \label{gensecexample}
Let $Q:\xymatrix@=5mm{0 \ar[r] & 1 & 2 \ar[l] & 3 \ar[l] & \ar[l] \cdots & \ar[l] n }$
and $c=(m,n,n-1,n-2,\cdots,1)$ for $m,n>0$. Let $Mat_{m\times n}$ be
the  space of $m\times n$ matrices, $B_{n}$ the group of invertible upper
triangular matrices and $\Gl_{m}\times B_{n}$ act on $Mat_{m\times n}$
by $(g,b)\cdot m = bmg^{-1}$. There is an inclusion of vector spaces 
$$Mat_{m\times n}\rightarrow \Rep(Q,c)$$ where a matrix 
$m\in Mat_{m\times n}$ is mapped to the representation with $m$ on 
the arrow $0\rightarrow 1$, and $\oplus_{i=1}^n Ae_i$ on the subquiver supported
on $1,2,3,\cdots, n$. The stabiliser of $\oplus_{i=1}^n Ae_i$ is $B_n$
and representations of $Q$ which restricts to  $\oplus_{i=1}^n Ae_i$
form an open subset of $\Rep(Q,c)$, so 
the $\Gl_{m}\times B_{n}$-space $Mat_{m\times n}$
is a generic section 
of the $\Gl(c)$-space $\Rep(Q,c)$.
\end{example}

\begin{example}
Let $Q:\xymatrix@=5mm{ 1\ar@/_/[r] \ar@/^/[r] & 2}$ be the Kronecker quiver 
and consider $\Rep (Q,c)$ with $c=(n, n)$. The $\Sigma_n$-space $k^n$ is a
generic section of the $\Gl(c)$-space $\Rep (Q, c)$, where $\Sigma_n$ is the 
symmetric group of order $n$ which acts on $k^n$ by the usual permutation of 
coordinates.
\end{example}

The generic section in the example, as well as the generic sections we consider
in this paper, preserve more information about the orbits than 
what we can deduce from Definition \ref{gensectiondef} alone. However,
this definition will be sufficient to prove the main result of this paper. 

\begin{lemma} \label{openlemma}
Let $V$ be a $G$-space and let the $H$-space $W$ be a generic section in $V$.
Then there is an open $G$-orbit in $V$ if and only if there is an open $H$-orbit
in $W$.
\end{lemma}
\begin{proof}
We may assume that the map $\phi:W\rightarrow V$ is an inclusion
$W\subseteq V$ of affine spaces. First, assume that $G\cdot x\subseteq V$
is open for some $x\in V$. By (2) of Definition \ref{gensectiondef},
$(G\cdot x)\cap W\subseteq W$ is nonempty and open. Then by (3), there
exists $x'\in (G\cdot x) \cap W'$ with $(G\cdot x')\cap W=H\cdot x'$.
This shows that $W$ has an open $H$-orbit.

Conversely, assume that $H\cdot w\subseteq W$ is open. We may assume
that $w\in W'$ and so $(G\cdot w) \cap W=H\cdot w$. Since the map 
$G\times W\ra V$ is dominant by Definition \ref{gensectiondef} part (2), 
the restriction $G\times (H\cdot w) \ra V$ is dominant, by comparison
of fibre dimensions. This shows that $V$ has an $G$-open orbit. 
\end{proof}

In particular, the lemma proves that $\{x\}\subseteq V$ with the trivial action, 
is a generic section of the $G$-space $V$ if and only if the orbit 
$G\cdot x\subseteq V$ is open. 

Generic sections define a relation on spaces with group actions, and two spaces
are said to be {\it generically equivalent} if they are related by a sequence of
generic sections. Two spaces with a dense open orbit are generically equivalent.

The key technical result in this paper, proved in Section 6, is that for certain 
quivers $Q$ and a projective representation $P$ of $Q$, the $\Aut_A P$-space
$\radE_A P$ is generically equivalent to a $\Gl(d')$-space $\Rep (Q',c')$, 
where $Q'$ is a quiver with the same underlying graph as $Q$, and $c'$ is a 
dimension vector constructed from $Q$ and $P$. We do not know if 
generic equivalences exist for all quivers $Q$, but they do exist for 
Dynkin quivers.  As a consequence, the main theorem 
stated in the introduction will follow by the Theorem \ref{Gabriel}, 
due to Gabriel. 

\section{Double quivers, adjoint actions and two equivariantly isomorphic varieties}

In this section we recall the construction of a  
finite dimensional quasi-hereditary algebra $D$, which can be used to
study $\Aut_A P$-orbits in $\radE_A P$ \cite{HV}.  We then recall and 
discuss some relevant properties of $D$ and its representation varieties.

Let $\tilde{Q}$ be the double quiver of $Q$, i.e. $\tilde{Q}_0=Q_0$ and 
$\tilde{Q}_1=Q_1\cup Q_1^*$ with $Q_1^*=\{\alpha^*: i\rightarrow j| 
\alpha: j\rightarrow i\in Q_1\}$. Let $\mathcal{I}$ be the ideal of $k\tilde{Q}$ 
generated by
$$\alpha^*\alpha - \sum_{\beta\in Q_1, t(\beta)=s(\alpha)} \beta\beta^*$$ 
for any arrow $\alpha\in Q_1$, and $$\alpha^*\beta$$ for pairs of arrows
$\alpha\neq \beta$ in $Q_1$ terminating at the same vertex. The algebra $D$ 
is defined as $$D=k\tilde{Q}/\mathcal{I}.$$

We now define a grading on $D$ which will be used in the next section. For a 
path $q$ in $Q$, we define $q^*=\alpha_n^*\dots \alpha_1^*$ if 
$q=\alpha_1\dots \alpha_n$ and $e_i=e_i^*$. The paths of the 
form $ab^*$ for paths $a,b$ in $Q$ with $s(a)=s(b)$ form a basis for $D$.
For any path $p$ in $\tilde{Q}$, let  $l^*(p)$ be the number of arrows from 
$Q_1^{*}$ in $p$. As the defining relations of $D$ are homogeneous with 
respect to $l^*$, we have the following lemma.

\begin{lemma}\label{grading}
$D$ is graded with respect to $l^*$.
\end{lemma}

Let $c=\dimv P$ be the dimension vector of $P$, $\Rep (D,c)$ the variety of 
$(\tilde{Q}, \mathcal{I})$-representations with dimension vector $c$ and let 
$\Gl(c)$ act on $\Rep (D, c)$ by conjugation. As $A=kQ$ is a subalgebra of $D$, 
any $D$-module can be considered as an $A$-module. We call a $D$-module 
that is  projective as an $A$-module an {\it $A$-projective $D$-module}. Let
$$\Rep (D,P)\subseteq \Rep (D,c)$$ be the subvariety consisting of the 
representations $X$ with $_{A}X=P$. The group $\Aut_A P$ acts on $\Rep(D,P)$
and orbits correspond to isomorphism classes of $A$-projective $D$-modules
$X$ with $_{A}X\cong P$.  

\begin{lemma} \label{irreducible}
\begin{itemize}
\item[]
\item[(i)] $\Rep (D,P)\subseteq \Rep (\tilde{Q},c)$ is an affine subspace.
\item[(ii)] $\Gl(c)\cdot \Rep (D,P)\subseteq \Rep (D,c)$ is irreducible
and open.
\end{itemize}
\end{lemma}
\begin{proof}
Note that the maps on each arrow in $Q$ are 
fixed for all point in $\Rep (D, P)$.  So $\Rep (D, P)$ is the solution space of the 
system of linear equations given by the defining relations of $D$. Therefore 
$\Rep (D, P)$ is an affine subspace of $\Rep(\tilde{Q},c)$. This proves (i)  and that $\Gl(c)\cdot \Rep (D,P)$ 
is irreducible.

Next observe that there is a natural morphism $\Rep (D,c)\rightarrow \Rep (Q,c)$ 
obtained by forgetting the $Q_1^*$-structure, i.e. the maps on each arrow in 
$Q_1^*$. Now $\Gl(c)\cdot \Rep (D,P)$ is the preimage of the open orbit of 
$P$, and so it is open. Thus (ii) follows.
\end{proof}

The following theorem can be deduced from its categorical version, Theorem 1.1 
in \cite{HV}. 
This result enables us to explore the existence of open $\Aut_A P$-orbits in 
$\radE_A P$ using representations of quivers.

\begin{theorem}\label{equivisom}
The varieties $\Rep (D, P)$ and $\radE_A P$ are $\Aut_A P$-equivariantly isomorphic.
\end{theorem}

Consequently, there is a one-to-one correspondence between $\Aut_A P$-orbits in
$\radE_A P$ and $(\tilde{Q}, \mathcal{I})$-representations that are isomorphic 
to $P$ as $A$-modules.

We now introduce two concepts to be used later. 
Let $P_i=Ae_i$ be the indecomposable projective $A$-module associated to 
the vertex $i\in Q_0$. Given an $A$-projective $D$-module $X$, denote 
by $d_i$ the multiplicity of $P_i$ as a summand in $_{A}X$ and let  
the {\it $\Delta$-dimension vector} of $X$ be defined as $$\dimv_\Delta X=(d_i)_i
\in  \mathbb{N}^n.$$ 
We let $$\Supp_{\Delta}X=\{i\in Q_0| d_i>0\}$$ 
and call it the {\it $\Delta$-support of $X$}. Note that the $\Delta$-support is 
in general different from the usual support $\Supp X$ defined using the dimension 
vector $\dimv X$. Given $d=(d_i)_i\in \mathbb{N}^n$,
let $$P(d)=\bigoplus_{i=1}^n{P_i}^{d_i}.$$

\begin{remark}
As mentioned earlier, the algebra $D$ is a quasi-hereditary algebra with 
the indecomposable projective $A$-modules $Ae_i$ as Verma modules \cite{HV}. 
The terminology $\Delta$-dimension vector and $\Delta$-support  
coincides with the one used 
in the setting of quasi-hereditary algebras, for example in \cite{BHRR, JSY}.
\end{remark}

\section{An exact sequence and a relative Voigt's lemma}

Let $J$ be the ideal in $D$ generated by the arrows in $Q^*_1$. Note that 
$D$ is a split extension of $A$ by $J$.  We construct two useful exact
sequences and give a criterion on the rigidity of  $A$-projective $D$-modules, 
analogues to Voigt's Lemma \cite{Voigt}.

\begin{lemma}\label{Jproj}
As  $D$-modules, $J\cong \oplus_{\beta\in Q^*, t(\beta)=i}De_i$
\end{lemma}

\begin{proof}
Observe that $J$ is the kernel of the natural surjective map
$$
D\rightarrow A.
$$ 
Now the lemma follows from the exact sequences (see \cite{HV, D}),
$$
\xymatrix{0 \ar[r] & \oplus_{\beta: i\rightarrow j \in Q^*}De_j\ar[r] &
De_{i} \ar[r] & Ae_i\ar[r] &0}, 
$$ 
for any vertex $i\in Q_0$.
\end{proof}

Recall the standard projective resolution of a representation $M$ of $Q$ 
(see \cite{CW}),
$$ 
\xymatrix{0 \ar[r] & \oplus_{\alpha\in Q_1} Ae_{t(\alpha)}
\otimes_k e_{s(\alpha)}M \ar[r] & \oplus_{i\in Q_0} Ae_{i}\otimes e_iM  \ar[r]&  
M\ar[r] &0.} 
$$
This resolution can be interpreted as a short exact sequence
$$
\xymatrix{0 \ar[r] &J'\otimes_SM \ar[r] & A\otimes_SM  \ar[r]&  M\ar[r] &0}, 
$$
where $S=\oplus_{i\in Q_0}ke_i$ and $J'$ is the kernel of the natural projection
$A\rightarrow S$. In other words, $J'$ is the Jacobson radical of $A$, which 
is the ideal of $A$ generated by arrows in $Q$. We have an analogous
resolution for $A$-projective $D$-modules as follows. 
The proof is similar to the proof of the standard resolution for quivers (see \cite{CW}).

\begin{lemma}\label{projresol}
Let $X$ be an $A$-projective $D$-module. Then the following exact
sequence is a projective resolution of $X$ as  a $D$-module,
$$\xymatrix{0\ar[r]& J\otimes_A X\ar[r]^f& D\otimes_A X \ar[r]^{\; \; \; \; \; \; g}
& X\ar[r]& 0 },$$ where $g(d\otimes m)=dm$ and $f(d\beta\otimes x)=
d\beta\otimes x-d \otimes \beta  x$ for $d\beta\in J, \beta\in
Q_1^{*}$, and $x\in X$.
\end{lemma}
\begin{proof}
Since $X$ is a projective $A$-module and $D$ and $J$ are projective 
$D$-modules, both $J\otimes_A X$ and $D\otimes_A X$ are projective
$D$-modules. So we only need to prove that the sequence is exact. By 
the definition of $f$ and $g$, the map $g$ is surjective and $gf=0$.
By applying $-\otimes_A X$ to the exact sequence
$$\xymatrix{0 \ar[r] & J \ar[r] & D \ar[r] &  A \ar[r] &0},$$
we obtain the exact sequence
$$
\xymatrix{0\ar[r]& J\otimes_A X\ar[r] & D\otimes_A X
\ar[r] & _AX\ar[r]& 0 },
$$
and so $\dim_k J\otimes_A X=\dim_k \Ker g$.  Therefore it remains to show 
that $f$ is injective. Let $\sum d\beta^*\otimes x$, for $\beta\in Q_1$ be in 
the kernel of $f$. Since $D$ is $l^*$-graded by Lemma \ref{grading},
we may assume that $l^*(d\beta^*)$ is constant on the terms 
in the sum. 
Now $f(\sum d\beta^*\otimes x)=\sum (d\beta^*\otimes x-d\otimes \beta^*x)=0$, which 
shows that $\sum d\beta^*\otimes x=0$, since 
$l^*(d)<l^*(d\beta^*)$. Hence $f$ is injective, and we are done.
\end{proof}

\begin{lemma}\label{variety1}
Let $X$ be an $A$-projective $D$-module. Then
\begin{itemize}
\item[(1)] $\Hom_D(D\otimes_AX, X)\cong \End_AX$.
\item[(2)] $\Hom_D(J\otimes _AX, X)\cong \radE_A X$.
\item[(3)] we have an exact sequence
$$
\xymatrix{0\ar[r]& \End_DX\ar[r]& \End_AX\ar[r]& \radE_AX\ar[r]&
\Ext^1_D(X,X) \ar[r]& 0. }
$$
\end{itemize}
\end{lemma}

\begin{proof}
First, $\Hom_D(D\otimes_AX, X)\cong \Hom_A(X, \Hom_D(D, X))\cong 
\Hom_A(X, X)=\End_AX.$ This proves (1).

Observe that for any indecomposable projective $A$-module $Ae_i$,
$J\otimes _A Ae_i\cong Je_i$. We have $J=\oplus_{\beta\in Q_1}D\beta^*$, and so $\Hom_D(J\otimes _AAe_ i, Ae_j)\cong \Hom_D(Je_i, Ae_j)\cong  \oplus_{s\rightarrow i}\Hom_D(De_s, Ae_j)
\cong\oplus_{s\rightarrow i}\Hom_A(Ae_s, Ae_j)\cong \radH_A(Ae_i, Ae_j)$,
where the sums are over all the arrows in $Q$ terminating at $i$. Therefore $$\Hom_D(J\otimes _AX, X)\cong \radE_A X.$$ This proves (2).

Applying $\Hom_D(-,X)$ to the sequence in Lemma \ref{projresol} 
gives the exact sequence,
$$
\xymatrix{0\ar[r]& \End_DX\ar[r]& \Hom_D(D\otimes X, X)
\ar[r]&\Hom_D(J\otimes X, X)   \ar[r]& \Ext^1_D(X,X)\ar[r]& 0. }
$$
Now by (1) and (2), we obtain the sequence in (3).
\end{proof}

Recall that a $D$-module is said to be {\it rigid} if $\Ext_D^1(X, X)=0$.
We can now use Lemma \ref{variety1} to prove the main result of this section on
the rigidity of $A$-projective $D$-modules, similar to Voigt's Lemma \cite{Voigt}. 
Moreover, it includes an inequality for the dimension of stabilizers, which is very useful
for determining whether an $A$-projective $D$-module is rigid or not.

\begin{theorem} \label{charlemma} 
Let $X$ be an $A$-projective $D$-module $X$ with $_{A}X=P(d)$ for a 
$\Delta$-dimension vector $d\in \mathbb{N}^n$. Then

\begin{itemize} \item[(1)] $\dim_k \End_DX\geq \sum_i
d_i^2$.
\item[(2)] the following are equivalent.
\begin{itemize}
\item[(i)] $Aut_A P\cdot X\subseteq \Rep (D,P)$ is open.
\item[(ii)] $X$ is rigid.
\item[(iii)] $\dim_k \End_DX=\sum_i d_i^2$.
\end{itemize}
\end{itemize}
\end{theorem}
\begin{proof}
First observe that  $\dim_k \End_AX-\dim_k \radE_AX=\sum {d_i}^2$. Now by
 Lemma \ref{variety1} $$\dim_k \End_DX= \dim_k \End_AX-\dim_k \radE_AX+\dim_k \Ext_D^1(X, X) \geq \sum_i d_i^2. $$ This proves (1). 
Moreover, the  equality holds  if and only if $\Ext_D^1(X,X)=0$. Thus the 
equivalence of (ii) and (iii) follows.

By Theorem \ref{equivisom}, $\dim_k \Rep (D, P)=\dim_k \radE_AX$. So
$$\dim_k Aut_A P\cdot X =\dim_k \End_AX-\dim_k \End_DX=\dim_k\Rep (D,P)-
\dim_k \Ext_D^1(X,X).$$  Therefore $i)$ and $ii)$ are equivalent.
\end{proof}

\section{Type $\mathbb{A}$}\label{typeA}

As preparation for the next section, we recall results on rigid $D$-modules 
in type $\mathbb{A}$ from \cite{BHRR, JSY}. At the same time we use 
Theorem \ref{charlemma} to give new, simpler and more uniform proofs 
of the two main results in \cite{BHRR} and \cite{JSY}, respectively. We 
also introduce a filtration on $D$, describe homomorphisms
in terms of this filtration and define  orders on indecomposable rigid $A$-projective 
$D$-modules, which will play a crucial role in the proofs of the main technical 
result in Section 6. Throughout this section $Q$ is a quiver of type 
$\mathbb{A}$, with vertices $\{1,\cdots, n\}$, and arrows 
$\alpha_i: i\rightarrow i+1$ or $\alpha_i:i\leftarrow i+1$ for $i=1,\cdots, n-1$.

\subsection{Linear $\mathbb{A}_n$} \label{linearA}

In this subsection $Q$ is linearly oriented and the arrows are  $\alpha_i:i
\rightarrow i+1$ for $i=1,\cdots,n-1$. In this case $\Aut_A P$ is a parabolic and 
the existence of a dense $\Aut_A P$-orbit in $\Rep(D,P)$ follows from the 
classical result of Richardson \cite{richardson}. We now recall the explicit 
construction of representations with dense open orbits by Br\"ustle, Hille, Ringel 
and R\"orhle \cite{BHRR}.

The projective $D$-module $Q_n=De_n$ at vertex $n$ is injective and has a
basis consisting of paths $pq^*$, where $q$ is  a path in $Q$ ending at the 
vertex $n$. A submodule $X$ of $Q_n$ has a basis given by a subset of the 
paths $pq^*$, and it is uniquely determined by its $A$-structure $_{A}X\cong
\oplus^n_{i=1}P_i^{d_i}$ with $d_i \in \{0,1\}$. Thus there is 
a natural bijection between subsets $I\subseteq Q_0$ and submodules of 
$Q_n$. More precisely, under this bijection a subset $I$ corresponds to the 
unique submodule of $Q_n$ with $\Delta$-support $I$. 
Let $X(I)$ denote the submodule of $Q_n$ corresponding to the subset $I$.
For any vector $d\in \mathbb{N}^n$, define $X(d)=\sum^t_{i=1} X(I_i)$, with $\dimv_\Delta X(d)=d$ and $I_1\subseteq I_2 \subseteq \cdots \subseteq I_t$.
We give an example to illustrate the construction.

\begin{example}\label{example2}
Let $n=3$ and $d=(2,1,2)$. The algebra $D$ is given by the quiver 
$$\xymatrix@=5mm{1 \ar@/^/[r]^{\alpha_1}
&2\ar@/^/[r]^{\alpha_2} \ar@/^/[l]^{\alpha_1^*}
& 3 \ar@/^/[l]^{\alpha_2^*} },$$ with
the ideal $\mathcal{I}$ generated by  $\alpha_1^*\alpha_1$ and
$\alpha_1\alpha_1^*-\alpha_2^*\alpha_2$.
The projective-injective D-module $Q_3$ has
the following $7$ nonzero submodules, where the first one is $Q_3$,
$$
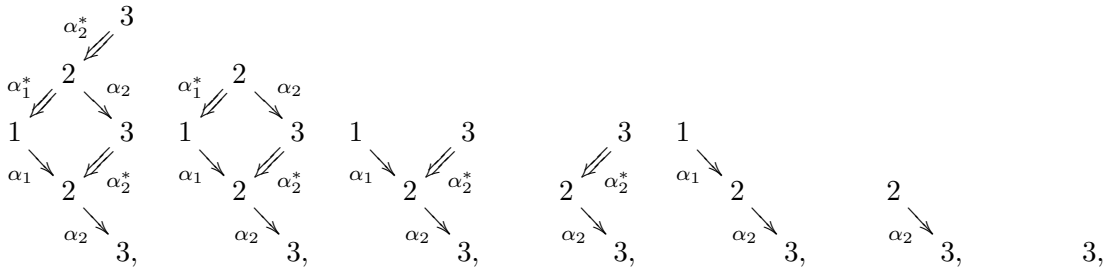
\xymatrix@=3mm{ && 3\ar@{=>}[dl]_{\alpha_2^*} \\
& 2\ar@{=>}[dl]_{\alpha_1^*}\ar[dr]^{\alpha_2} & & & 2\ar@{=>}[dl]_{\alpha_1^*}\ar[dr]^{\alpha_2} \\ 1\ar[dr]_{\alpha_1}
&& 3\ar@{=>}[dl]^{\alpha_2^*} & 1\ar[dr]_{\alpha_1}  
&& 3\ar@{=>}[dl]^{\alpha_2^*} & 1\ar[dr]_{\alpha_1}  
&& 3\ar@{=>}[dl]^{\alpha_2^*}
&&& 3\ar@{=>}[dl]^{\alpha_2^*} & 1\ar[dr]_{\alpha_1}\\
& 2\ar[dr]_{\alpha_2} & && 2\ar[dr]_{\alpha_2}  &&& 2\ar[dr]_{\alpha_2}  
&&&2\ar[dr]_{\alpha_2} &&& 2\ar[dr]_{\alpha_2}  &&& 
2\ar[dr]_{\alpha_2} \\ && 3,&&&3,&&&3,&&&3,&&&3,&&&3, &&& 3,
}$$
corresponding to the subsets $\{1, 2,3\}$, $\{1, 2\}$, $\{1, 3\}$, $\{2, 3\}$, 
$\{1\}$, $\{2\} $, $\{3\}$, respectively. In the picture a number $i$ indicate
a one dimensional basis element at vertex $i$ 
and the arrows indicate the nonzero action of the arrows in $\tilde{Q}_{1}$. 

For $d=(2, 1, 2)$, we have $X(d)=X(I_1)\oplus X(I_2)$ for the subsets
$I_1=\{1, 2, 3\}$ and $I_2=\{1, 3\}$, respectively.
\end{example}

\begin{theorem}[\cite{BHRR}] \label{BHRR}
An $A$-projective $D$-module $X$ is rigid if and only if $X\cong X(d)$
for some $d\in \N^{n}$.
\end{theorem}

We will give a new proof of the above theorem. First recall a 
lemma in \cite{BHRR} on the dimension of homomorphism spaces between
submodules of $Q_n$, and include a new proof. 

\begin{lemma}[\cite{BHRR}]\label{dim1}
Let $X(I)$ and $X(J)$ be submodules of $Q_n$ with $\Delta$-support $I$ and $J$, 
respectively, and $I\subseteq J$. Then $$\dim_k \Hom_D(X(I), X(J))=\dim_k 
\Hom_D(X(J),X(I))= |I|.$$
\end{lemma}
\begin{proof}
Note that the module $Q_n$ is generated by $e_n$ and $\End_D(Q_n)$ is an 
$n$-dimensional vector space with basis $f_0, \dots, f_{n-1}$, where
$f_i(e_n)=(\alpha_{n-1} \alpha_{n-1}^*)^{i}e_n$ for $0\leq i\leq n-1$.
In particular, $f_0$ is the identity map.
For any two submodules $X$ and $Y$ of $Q_n$ and any homomorphism 
$f:X \ra Y$, since $Q_n$ is injective,
we have the following commutative diagram for a $g\in \End_DQ_n$
$$
\xymatrix{
X \ar[r]^{\subseteq} \ar[d]_f& Q_n\ar@{.>}[d]^g
\\
Y\ar[r]^{\subseteq} &Q_n.}
$$
That is, $f=g|_{X(I)}$ and so
$f$ is a linear combination of the restrictions of $f_0, \dots, f_{n-1}$ to $X(I)$.

By the definition, the nonzero restrictions of $f_0, \dots, f_{n-1}$ in
$\Hom_D(X(I), X(J))$ are linearly independent.
Moreover, $f_i|_{X(I)}=0$ if and only if $i\geq \sum d_i$,
where $(d_i)_i=\dimv_\Delta X$.
Therefore $$\dim_k \Hom_D(X(I),X(J)) = \sum_i d_i = |I|.$$

Similarly,  $\dim_k\Hom_D(X(J), X(I))=|I|$.
\end{proof}

\begin{proof}[Proof of Theorem \ref{BHRR}]
By Theorem \ref{charlemma} and Lemma \ref{dim1}, each submodule $X(I)$ is 
rigid. So by the construction of $X(d)$, to show that $X(d)$ is rigid, we need 
only prove that $$X=X(I)\oplus X(J)$$ with $I\subseteq J$ is rigid. We have
$_{A}X\cong \oplus_iP(i)^{d_i}$, where $d_i=2$ if $i\in I\cap J$, $d_i=1$ if
$i\in J \backslash I$ and $d_i=0$ if $i\not\in J$. By Lemma \ref{dim1},
$\dim_k \End_D(X)=3|I|+|J|=4|I|+|J\backslash  I|=\sum d_i^2$, and so $X$ is 
rigid by Theorem \ref{charlemma}. This proves
the existence of an open orbit in $\Rep (D, P)$ for any $P$. Now $\Rep (D, P)$ 
is irreducible by Lemma \ref{irreducible}, and so by 
Theorem \ref{charlemma}, any rigid $D$-module is isomorphic to $X(d)$ 
for some $d\in \N^n$.
\end{proof}

\subsection{Gluing modules at sinks and sources}

In this subsection $Q$ is of type $\mathbb{A}$ with arbitrary
orientation.
We recall how to glue a pair of $A$-projective $D$-modules at an admissible interior
vertex to obtain a new $A$-projective $D$-module \cite{JSY}.

Let $$i_1<i_2<\cdots <i_{t-1}<i_t$$ be the complete list of interior admissible vertices 
in $Q$ and let $i_0=1$ and $i_{t+1}=n$. Let $u=i_l$ 
be one of the interior admissible vertices, and let $M'$ and $M''$ be two indecomposable
$A$-projective $D$-modules with $\Supp_{\Delta}(M')\subseteq$ $\{1,\dots,u\}$, 
$\Supp_{\Delta}(M'')\subseteq$ $\{u,\dots,n\}$ and $(\dimv_{\Delta} M')_u=1=
(\dimv_{\Delta} M'')_u$.

Assume first that $u$ is a sink in $Q$. Let $N'$ be the $D$-submodule of $M'$ 
generated by $M'_j$ for all $j<u$. Then $N'$ is not $\Delta$-supported at $u$, 
i.e. $(\dimv_{\Delta}N')_u=0$. Since $u$ is a sink and $(\dimv_{\Delta} M')_u=1$, 
we have the short exact sequence
$$\xymatrix{0 \ar[r] &N' \ar[r] & M' \ar[r]^{f'} & P_u \ar[r] & 0,} $$ where 
$P_u=S_u$ is the simple projective $A$-module associated to vertex $u$. 
Similarly,  $$\xymatrix{0 \ar[r] & N''\ar[r] & M'' \ar[r]^{f''} & P_u \ar[r] & 0,}$$
where $N''\subseteq M''$ is the submodule generated by $M''_j$ for all $j>u$.
Let $M$ be given by the pullback of $f'$ and $f''$, that is, we have a short exact 
sequence $$\xymatrix{0 \ar[r] & M \ar[r] & M'\oplus M'' \ar[r]
& P_u \ar[r] & 0.}$$  We say that $M$ is obtained by {\it gluing} $M'$ and $M''$ at 
the sink $u$.

We now define gluing of homomorphisms at the sink $u$. Assume that $K$ is 
obtained by gluing $K'$ and $K''$ at $u$ and let $L'\subseteq K'$ and 
$L''\subseteq K''$ be the submodules generated by the spaces $K'_i$ for $i<u$ 
and $K''_i$ for $i>u$, respectively.  Let $g':M'\ra K'$ and $g'':M''\ra K''$ be 
homomorphisms of $D$-modules.
Then $g'(N')\subseteq L'$ and $g''(N'')\subseteq L''$, and so there are induced 
maps $\tilde{g'},\tilde{g''}:P_u\ra P_u$. Using the pullback sequence defining 
$M$, we see that there is an induced map $g:M\ra K$ if $\tilde{g'}=\tilde{g''}$. 
In this case we say that $g$ is obtained by gluing $g'$ and $g''$ at the sink $u$. 

Now assume that $u$ is a source. Note that in this case, the projective $A$-module 
$P_u$ is also projective when viewed as a $D$-module. We have short exact sequences 
$$\xymatrix{0 \ar[r] & P_u\ar[r]^{f''} & M' \ar[r] & N' \ar[r] & 0,}$$
$$\xymatrix{0 \ar[r] & P_u\ar[r]^{f''} & M'' \ar[r] & N'' \ar[r] & 0,}$$
where $P_u\subseteq M'$ and $P_u\subseteq M''$ are the submodules generated 
by $M'_u$ and $M''_u$, respectively. Let $M$ be given by the pushout of $f'$ 
and $f''$, 
$$\xymatrix{0 \ar[r] & P_u \ar[r] & M'\oplus M'' \ar[r] & M \ar[r] & 0.}$$ 
We then say that $M$ is obtained by gluing $M'$ and $M''$ at the source $u$.
Assume $K$ is obtained by gluing $K'$ and $K''$ at the source $u$.  Given a pair 
of homomorphisms $g':M'\ra K'$ and $g'':M''\ra K''$ there is an induced map $g:M\ra K$ if 
the induced maps $\tilde{g'}=g'|_{P_u}$ and $\tilde{g''}=g''|_{P_u}$ are equal. 
We then say that $g$ is obtained by gluing $g'$ and $g''$ at the source $u$.

\begin{lemma} \label{gluingmap}
Let $u$ be an interior admissible vertex. Then
$$\Hom_D(M,K)\cong \{(g',g'')\in \Hom_D(M',K')\oplus \Hom_D(M'',K'') |\tilde{g'}=
\tilde{g''}\}$$ 
\end{lemma}
\begin{proof}
Assume $u$ is a sink. Let $$g:M\ra K$$ be a homomorphism. Then there are 
the short exact sequences 
$$\xymatrix{0 \ar[r] & N' \ar[r] & M \ar[r] & M'' \ar[r] & 0}$$ and
$$\xymatrix{0 \ar[r] & L' \ar[r] & K \ar[r] & K'' \ar[r] & 0,}$$ where
$N'\subseteq M$ and $L'\subseteq K$ are the submodule generated by 
$M_i$ and $K_i$, respectively, for $i<u$.  Then $g(N')\subseteq L'$ and
so there is an induced map $g'':M''\ra K''$ and a commutative diagram
$$\xymatrix{0 \ar[r] & N' \ar[r] \ar[d]^{g|_{N'}} & M \ar[r] \ar[d]^g & M'' 
\ar[r] \ar[d]^{g''}&  0 \\
0 \ar[r] & L' \ar[r] & K \ar[r] & K'' \ar[r] & 0.}$$
Similarly, there is a map $g':M'\ra K'$. Any map $g:M\ra K$ is therefore
obtained by gluing maps $g'$ and $g''$ at $u$, and so there is an injection of 
vector spaces $\Hom_D(M,K)\ra \Hom_D(M',K')\oplus \Hom_D(M'',K'')$ given by
$$g\mapsto (g',g'').$$ The image is equal to the pairs of maps $(g',g'')$
with $\tilde{g'}=\tilde{g''}$, since each such pair can be glued at $u$.

This proves the lemma in the case where $u$ is a sink. The case where 
$u$ is a source is similar and is left to the reader.
\end{proof}

\subsection{The construction of rigid modules}
For $d\in \N^n$ and any $0\leq s\leq t$, let 
$d^s\in \N^n$ be the vector defined by $(d^s)_j=d_j$ for $j=i_{s},i_{s}+1,\cdots,i_{s+1}$ and $(d^{s})_j=0$, 
otherwise. Using Theorem \ref{BHRR}, we can construct a rigid representation 
with $\Delta$-dimension vector $d^s$ for the double quiver supported on 
$\{i_s,\cdots,i_{s+1}\}$ subject to the corresponding relations. This module,
which we denote by $Y(d^s)$, may not be $A$-projective when considered
as a $D$-module. However, if $i_{s}$ is an interior source, we extend 
$Y(d^s)$ by $P_{i_{s}-1}^{d_{i_s}}$, and if $i_{s+1}$ is an interior 
source, we extend $Y(d^s)$ by $P_{i_{s+1}+1}^{d_{i_{s+1}}}$, to 
obtain an $A$-projective $D$-module, which we denote by $X(d^s)$. 
Since the extension preserves the $\Delta$-dimension vector and the 
dimension of the endomorphism ring, Theorem \ref{charlemma} shows that 
the $D$-module $X(d^s)$ is rigid. We will glue the indecomposable 
summands of $X(d^s)$ for all $s$ with respect to 
an order defined below to obtain a rigid module with $\Delta$-dimension 
vector $d$. 

Following the construction of $X(d^s)$ each indecomposable 
summand of $X(d^s)$ is completely 
determined up to isomorphism by its $\Delta$-support. Suppose that $M$ is 
an  indecomposable module obtained by gluing modules 
$X(I_s)$ at interior admissible vertices, where $I_s=\Supp_\Delta X\cap \{i_s,\dots, i_{s+1}\}$ for all $s$,
and that $N$ is an indecomposable module obtained by gluing $X(J_s)$,
where $J_s=\Supp_\Delta Y\cap \{i_s,\dots, i_{s+1}\}$ for all $s$.
For such indecomposable modules we define an order $\leq_u$ for
any vertex $u\in Q_0$.

\begin{definition}\label{order}
 Let $u$ be a vertex with $i_{v}<u\leq i_{v+1}$. 
Suppose that both $M$ and $N$ are supported 
(but not necessarily $\Delta$-supported) 
at $u$. We define $M\leq_u N$ if, for any $s$ with $I_s$ and $ J_s$ nonempty, 
$I_s\subseteq J_s$ if $s-v$ is even and $I_s\supseteq J_s$ if $s-v$ is odd. 
\end{definition}

The construction of a rigid $A$-projective module $X(d)$ with 
$\Delta$-dimension vector $d$ is done by induction on the number of interior 
admissible vertices in $Q$. Clearly, if there are no interior admissible vertices, 
then $Q$ is linearly oriented, and we are done.

Now suppose that the number of interior admissible vertices is $t>0$. For any 
$s\geq 1$, define a vector  $e^s$ given by $(e^s)_j=d_j$ if $j\leq i_s$ and
$(e^s)_j=0$ otherwise. We have $e^1=d^0$ and $e^{t+1}=d$. 
We suppose that $$X(d^0)=X^{01}\oplus \dots
\oplus X^{0d_{i_1}}\oplus X'$$ and $$X(d^1)=X^{11}\oplus \dots \oplus X^{1d_{i_1}}
\oplus X'',$$ where $i_1 \in \Supp_{\Delta}X^{01}\subseteq \dots\subseteq 
\Supp_{\Delta}X^{0d_{i_1}}$ and $\Supp_{\Delta}X^{11}\supseteq 
\dots \supseteq \Supp_{\Delta}X^{1d_{i_1}} \ni i_1$, each $X^{ij}$ is 
indecomposable and $X'$ and $X''$ are not $\Delta$-supported at $i_1$. 
For $1\leq s \leq d_{i_1}$, we glue $X^{0i}$ and $X^{1i}$ at $i_1$ and 
obtain a module we denote by $X^i$. Let $X(e^2)=X^1\oplus \dots \oplus X^{d_{i_1}}
\oplus X'\oplus X''$. By construction, the summands of $X(e^{2})$ that 
are $\Delta$-supported at $i_2$ are totally ordered with respect to 
$\leq_{i_2}$. 

By induction we may assume that we have constructed an $A$-projective 
$D$-module $X(e^s)$ with $\Delta$-dimension vector $e^s$ and that 
the indecomposable summands of $X(e^s)$ that have $\Delta$-support 
at $i_s$ are totally ordered using $\leq_{i_s}$.
Here, as in the case of $s=2$, by abuse of notation we again use 
$X^1\leq_{i_s}\dots \leq_{i_s} X^{d_{i_s}}$ to 
denote the summands of $X(e^s)$ that are $\Delta$-supported at $i_s$. 
We have $X(e^s)=X^1 \oplus \dots \oplus X^{d_{i_s}}\oplus Z$ and 
$X(d^{s})=X^{s1}\oplus \dots X^{sd_{i_s}}\oplus Z' $, where 
$\Supp_{\Delta} X^{s1}\supseteq \dots \supseteq \Supp_{\Delta}X^{sd_{i_s}}\ni i_{s}$.
For each $i$, glue $X^{i}$ with $X^{si}$ at vertex $i_s$ and obtain a 
new module, which we again denote by $X^i$. Now let $X(e^{s+1})=X^1
\oplus \dots\oplus X^{d_{i_s}}\oplus Z\oplus Z'$. 

\begin{lemma}
Suppose that $X(e^{s})$  above is rigid. Let $M$, $N$ be two 
summands of $X(e^{s})$, which are $\Delta$-supported at $i_s$. 
Assume that $M\leq_{i_s} N$ and 
$i_s\in I\subseteq J \subseteq \{i_{s},\cdots,i_{s+1}\}$. Let $X$ be
obtained by gluing $M$ with $X(J)$ and let $Y$ be obtained by gluing
$N$ with $X(I)$. Then $X\oplus Y$ is rigid.
\end{lemma}
\begin{proof}
We claim that $$\dim_k\Hom_D(X,Y)=\dim_k\Hom_D(M,N)+|I|-1.$$ We first consider 
the case where $i_s$ is a source. By Lemma \ref{gluingmap}, 
$$\dim_k\Hom_D(X,Y)\leq \dim_k\Hom_D(M,N)+ \dim_k\Hom_D(X(J),X(I)). $$
Again by Lemma \ref{gluingmap} any pair $$(g',g'')\in \Hom_D(M, N)\oplus H
om(X(J),X(I))$$ with $g'|_{M_{i_s}}=0$ and $g''|_{X(J)_{i_s}}=0$ can be 
glued to a homomorphism $X\ra Y$. Since $M\leq_{i_s} N$, there exists a map 
$g'_1:M\lra N$ with $g'_1|_{M_{i_s}}\neq 0$. Moreover, the subspace in 
$\Hom_D(M,N)$ consisting of the maps that are zero at $i_s$ has codimension $1$.
In $\Hom_D(X(J),X(I))$, the subspace consisting of the maps that are zero at $i_s$ 
has codimension at most $1$, with equality if and only if $I=J$. So 
$$\dim_k\Hom_D(X,Y)\geq  (\dim_k\Hom_D(M,N)-1)+ (\dim_k\Hom_D(X(J),X(I)) -1). $$

If $I=J$, then $g'_1:M\ra N$ can be glued to an isomorphism 
$g''_1:X(J)\ra X(I)$. Therefore 
$$\dim_k\Hom_D(X,Y)=  \dim_k\Hom_D(M,N) + \dim_k\Hom_D(X(J),X(I)-1. $$

If $I\subsetneq J$, then $g''|_{X(J)_{i_s}}=0$ for all maps 
$g''\in  \Hom_D( X(J), X(I))$. Thus $$\dim_k\Hom_D(X,Y) =  (\dim_k\Hom_D(M,N)-1)+ 
\dim_k\Hom_D(X(J),X(I)).$$ 

By Lemma \ref{dim1}, $\dim_k\Hom_D(X(J),X(I)) =|I|$ and so the claim follows for 
$i_s$ a source. The proof of the claim is similar for the case where $i_s$ is a 
sink and we leave it to the reader. 

Similarly, we have 
$$\dim_k\Hom_D(X,X)=\dim_k\Hom_D(M,M)+|I|-1,$$ 
$$\dim_k\Hom_D(Y,X)=\dim_k\Hom_D(N, M)+|I|-1,$$
$$\dim_k\Hom_D(Y,Y)=\dim_k\Hom_D(N,N)+|J|-1.$$

We show that $X\oplus Y$ is rigid. Suppose that 
$\dimv_{\Delta}(X\oplus Y)=(d_i)_i$. Then $d_{i}=2$ for 
$i\in I$, $d_j=1$ for $j\in J\backslash I$ and $d_i = 
(\dimv_\Delta M\oplus N)_i$ otherwise. 
 By our assumption, 
$M\oplus N$ is rigid and so by Theorem \ref{charlemma}, 
$$\dim_k \End_D(M\oplus N)=\sum_{i=1}^{i_s}d_i^2.$$ Therefore 
$$\dim_k \End_D(X\oplus Y)=\dim_k \End_D(M\oplus N)+ 4|I|+|J\backslash I|-4=
\sum_id_i^2.$$ Again by Theorem \ref{charlemma}, $X\oplus Y$ 
is rigid.
\end{proof}

Using the lemma we can finish the case of type $\mathbb{A}$.

\begin{theorem}[\cite{JSY}] \label{JSY} 
An $A$-projective $D$-module $X$ is rigid if and only if $X\cong X(d)$ for
some $d\in \mathbb{N}^n$.
\end{theorem}
\begin{proof}
We use induction on $s$ to prove that $X(e^{s+1})$ is rigid for any $s$, thus
so is $X(d)=X(e^{t+1})$. Note that $X(e^1)=X(d^0)$ is rigid, and 
assume that $X(e^s)$ is rigid. By the previous lemma we need
only to prove the rigidity of $X\oplus Y \oplus L \oplus L'$, where $X$ and $Y$ are
as in the previous lemma, $L$ and $L'$ are, respectively, summands of $X(e^s)$ and 
$X(d^s)$ without $\Delta$-support at $i_s$. Then $\Hom_D(L,X\oplus
Y)=\Hom_D(L,M\oplus N)$ and $\Hom_D(X\oplus Y,L)=\Hom_D(M\oplus N,L)$ and so
$L\oplus X\oplus Y$ is rigid, by Theorem \ref{charlemma}. Similarly,
for $L\oplus L'$ and $L'\oplus X\oplus Y$. Therefore $X\oplus Y \oplus L \oplus L'$ is 
rigid. 

Conversely, if $X\in \Rep(D,P(d))$ is rigid, then by Lemma \ref{irreducible} 
and Theorem \ref{charlemma}, $X$ is isomorphic to $X(d)$.
\end{proof}

We give an example of the above construction. 

\begin{example}\label{example3}
Let $Q$ be the quiver 
$$\xymatrix@=5mm{ 1 \ar[r]&2\ar[r] &3& 4\ar[l]&5\ar[l]\ar[r]&6\ar[r]&7.}$$

(1) Let $d=(0, 0, 2, 1, 2, 1, 2)$. Then $d^0=(0, 0, 2, 0, 0, 0  0)$, $d^1=(0, 0, 2, 1, 2, 0, 0)$ 
and  $d^2=(0, 0, 0, 0, 2, 1, 2)$, and $X(d^1)=M^1\oplus M^2$ and 
$X(d^2)=N^1\oplus N^2$ as follows.

$$
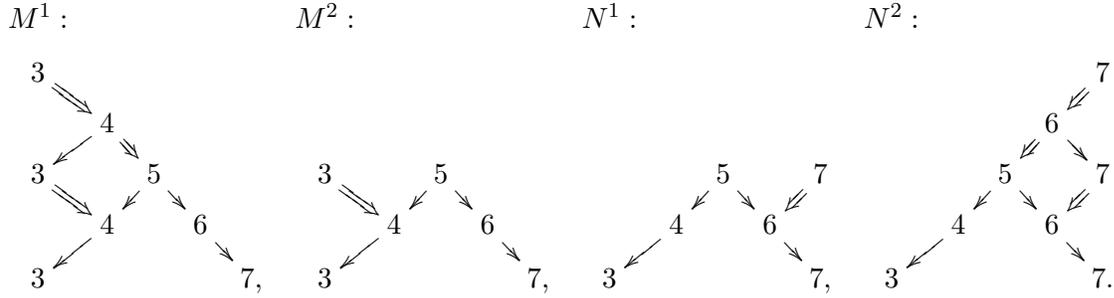
\xymatrix@=2mm{  
M^1: &&&&&  M^2:  &&&&&   N^1: &&&&&N^2:\\
3\ar@{=>}[dr] &&&&& & &&&  &&  & &&  &&&&& 7\ar@{=>}[dl] \\
& 4\ar[dl]\ar@{=>}[dr] &&&  &&  & &&  &&&&& &&&&6\ar[dr]\ar@{=>}[dl]\\
3\ar@{=>}[dr]&&5\ar[dl]\ar[dr]&&& 3\ar@{=>}[dr] && 5\ar[dr]\ar[dl] &&&& 
&5\ar[dl]\ar[dr] && 7\ar@{=>}[dl] & &&  5\ar[dl]\ar[dr]&&7\ar@{=>}[dl] \\
& 4\ar[dl]&&6\ar[dr]&&&4 \ar[dl]&& 6 \ar[dr]
&&& 4\ar[dl] && 6\ar[dr] && & 4\ar[dl]&&6\ar[dr] \\
3 &&&& 7, &3&&&&7,& 3&&&&7, &3 &&&&7. }$$

We have $M^1\geq_{5} M^2$, $N^1=X(\{5, 7\})$ and $N^2=X(\{5,6,7\})$. 
So $X(d)$ is the direct sum of the gluings of $M^1$, $M^2$ with $N^1$ and $N^2$, 
respectively as follows, 

$$
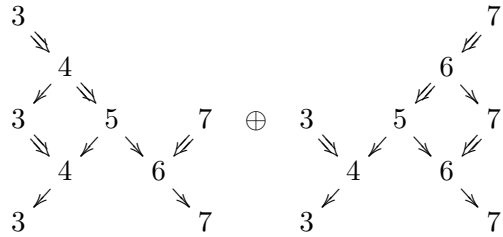
\xymatrix@=2mm{
3\ar@{=>}[dr] &&&&& & &&& &7\ar@{=>}[dl]
\\
& 4\ar[dl]\ar@{=>}[dr] 
&&&  && && &6\ar@{=>}[dl]\ar[dr]\\
3\ar@{=>}[dr]&&5\ar[dl]\ar[dr]&& 7\ar@{=>}[dl]& \oplus &
3\ar@{=>}[dr]&&5\ar[dr]\ar[dl]&&7\ar@{=>}[dl]
\\
&4\ar[dl] &&6\ar[dr]    &&&& 4\ar[dl]&& 6\ar[dr]\\
3 &&&& 7 &         &3 &&&&7
}
$$

(2) Let $d=(1,2,2,0,2,1,2)$. In this case $e^2=(1,2,2,0,2, 0, 0)$ and $d^2$ is the 
same as in (1). We have $X(e^2)=M\oplus N$ as follows. 
$$
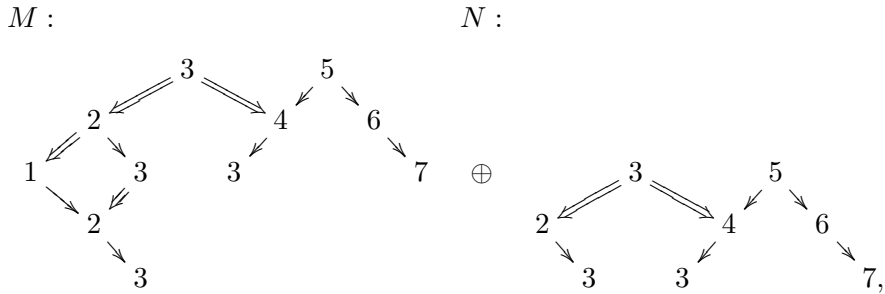
\xymatrix@=2mm{
M:  &&&&&& &&&N:\\
&&& 3\ar@{=>}[dll]\ar@{=>}[drr] &&& 5\ar[dl]\ar[dr]\\
& 2\ar@{=>}[dl]\ar[dr]&&&&4\ar[dl]&&6 \ar[dr]\\
1\ar[dr]&&3\ar@{=>}[dl] &&3&&& &7& \oplus &&&3 \ar@{=>}[dll]\ar@{=>}[drr]
&&&5\ar[dr]\ar[dl]\\
&2\ar[dr]&&&& && && &2\ar[dr]&&&&4\ar[dl] && 6\ar[dr]&&\\
&&3 &&& &&&&&&3&&3 &&&&7,}
$$
where $M\leq_5 N$. So $X(d)$ is the direct sum of the gluing of $M$ and $N$ with 
$N^2$ and $N^1$ from (1), respectively, as follows.

$$
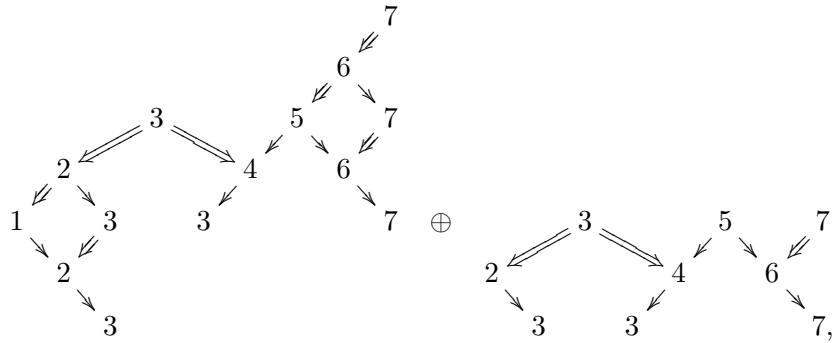
\xymatrix@=2mm{
&&&&&&& & 7\ar@{=>}[dl]\\
&&&&&&& 6\ar@{=>}[dl]\ar[dr]& 
\\
&&& 3\ar@{=>}[dll]\ar@{=>}[drr] &&& 5\ar[dl]\ar[dr] &&7\ar@{=>}[dl]\\
& 2\ar@{=>}[dl]\ar[dr]&&&&4\ar[dl]&&6 \ar[dr]\\
1\ar[dr]&&3\ar@{=>}[dl] &&3&&& &7& \oplus &&&3 \ar@{=>}[dll]\ar@{=>}[drr]
&&&5\ar[dr]\ar[dl] && 7\ar@{=>}[dl]\\
&2\ar[dr]&&&&& && & &2\ar[dr]&&&&4\ar[dl] && 6\ar[dr]&&\\
&&3 && &&&&&&&3&&3 &&&&7,}
$$
\end{example}

\subsection{A filtration on rigid modules} \label{grading}

We now discuss a filtration on $D$ and on $X(d)$ to be used in the 
next section.  Let $D^i\subseteq D$
be the subspace spanned by paths of the form   
$$a(\alpha\alpha^*)^ib^*$$
where $\alpha\in Q_1$ and $a$ and $b$ are paths in $Q$ with no arrow in common. Then $$D=\bigoplus_{i\geq 0}D^i$$ 
as vector spaces, and we say that the paths in $D^i$ have degree $i$. Let $$\mathcal{J}=\oplus_{i\geq 1}D^i,$$ be the 
ideal generated by all paths of the form $\alpha\alpha^*$ for $\alpha\in Q_1$. 
We have a filtration $$D=\mathcal{J}^0\supseteq \mathcal{J}
\supseteq \mathcal{J}^2 \supseteq \dots ,$$ such that
the inclusion $\mathcal{J}^i\subseteq D$ induces 
an isomorphism $$\mathcal{J}^i/\mathcal{J}^{i+1}\cong D^i.$$
The decomposition of $D$ induces a decomposition on any indecomposable rigid
$A$-projective $D$-module $$X=\oplus_{i\geq 0}X^i,$$ in such a
way that the component $X^i$ is identified with
$\mathcal{J}^iX/\mathcal{J}^{i+1}X$, in particular, $$X^0\cong
X/\mathcal{J}X.$$

\begin{example} The decomposition on the first summand of $X(d)$ in 
Example \ref{example3}(2) is as follows.
$$
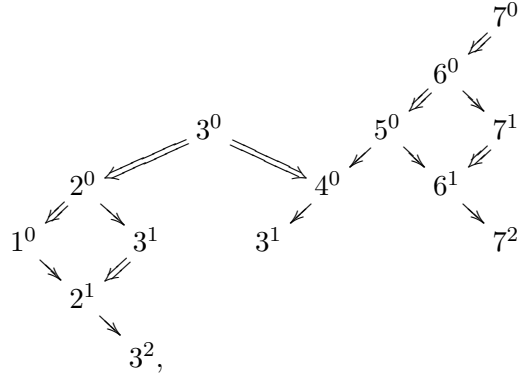
\xymatrix@=2mm{
&&&&&&& & 7^0\ar@{=>}[dl]\\
&&&&&&& 6^0\ar@{=>}[dl]\ar[dr]& \\
&&& 3^0\ar@{=>}[dll]\ar@{=>}[drr] &&& 5^0\ar[dl]\ar[dr] &&7^1
\ar@{=>}[dl]\\
& 2^0\ar@{=>}[dl]\ar[dr]&&&&4^0\ar[dl]&&6^1 \ar[dr]\\
1^0\ar[dr]&&3^1\ar@{=>}[dl] &&3^1&&& &7^2\\
&2^1\ar[dr]&&&&& && & \\
&&3^2, && &&&&&&}$$
where $i^j$ indicates a basis element of degree $j$ at vertex $i$.
\end{example}

We collect some basic properties in the following lemmas.

\begin{lemma} \label{gradingtool} 
Let $X$ be a rigid indecomposable $A$-projective $D$-module supported at $u$. 
We have $X^0_u\cong k$. Moreover,
\begin{itemize}
\item[(1)] if $u$ is a source, then $X^0_u=X_u$,
\item[(2)] if $u$ is an interior sink, then 
$X^i_u=(\alpha_{u-1}\alpha_{u-1}^*)^iX^0_u
\oplus (\alpha_{u}\alpha_{u}^*)^iX^0_u$ for $i\geq 1$, and
\item[(3)] if there is precisely one arrow $\alpha\in Q_1$
terminating at $u$, then $X^i_u=(\alpha\alpha^*)^iX^0_u$
\end{itemize}
\end{lemma}

Let $u$ be an interior vertex of $Q$ with $i_v<u\leq i_{v+1}$ and assume
that $i_v$ is a source and $i_{v+1}$ is a sink. Let
$M$ and $N$ be two indecomposable rigid $D$-modules supported (but
not necessarily $\Delta$-supported at $u$) such that $M<_u N$.
Let 
$$\Hom(M,N)^0_u=\{f|_{M^0_u}\mid f\in \Hom_D(M,N), 
f(M^0_u)\subseteq N^0_u\}$$
denote the space of maps obtained by restricting a map $f$ to $M^0_u$
where $f(M^0_u)\subseteq N^0_u$. 
We have $M^0_u=k=N^0_u$, and so both $\Hom(N,M)^0_u$ 
and $\Hom(M,N)^0_u$ are at most one dimensional.

\begin{lemma}\label{gradedhom}
Let $M$ and $N$ be as above. Then
\begin{itemize}
\item[(a)] $\Hom(N,M)^0_u=k$ if and only if $(\dimv_\Delta M)_w =
(\dimv_\Delta N)_w$ for all $w>u$.
\item[(b)] $\Hom(M,N)^0_u=k$ if and only if $(\dimv_\Delta M)_w =
(\dimv_\Delta N)_w$ for all $w\leq u$.
\end{itemize}
\end{lemma}

\begin{proof}
(a) Suppose that $(\dimv_\Delta M)_w=(\dimv_\Delta N)_w$ for all 
$w>u$. Let $I$ and $J$ be the  $\Delta$-supports of $M$ and $N$, 
respectively, and let  $I_s=I\cap [i_s, i_{s+1}]$ and  
$J_s=J\cap [i_s, i_{s+1}]$, where $[i_s, i_{s+1}]=\{i_s, i_s+1, 
\dots, i_{s+1}\}$. Since $M<_u N$ we have by the definition of $\leq_u$ 
that $I_s\subseteq J_s$ if $s-v$ is even and $I_s\supseteq J_s$ if $s-v$ is 
odd, whenever both $I_s$ and $J_s$ are nonempty. 
Since $(\dimv_\Delta M)_w =(\dimv_\Delta N)_w$ for all $w>u$, 
we have $I_s=J_s$ for $s\geq v+1$. 
By the construction of rigid modules in Theorem \ref{JSY}, both $M$ and 
$N$ are obtained by gluing indecomposable rigid $D$-modules $X(I_s)$ 
and $X(J_s)$, respectively. Let $t$ be the largest integer such that $I_t\neq J_t$.
We may glue identity maps $X(J_s)\rightarrow X(I_s)$ for $s>t$ and
a surjection $X(J_t)\rightarrow X(I_t)$ if $v-t$ is even, and an injection 
$X(J_t)\rightarrow X(I_t)$ if $v-t$ is odd, to obtain a homomorphism 
$f:N\rightarrow M$. In both cases $0\neq f(N^0_u)\subseteq M^0_u$, 
and so $\Hom(N,M)^0_u=k$.

Conversely, assume there is a homomorphism $f:N\rightarrow M$ with
$0\neq f(N^0_u)\subseteq M^0_u$. Since $i_{v+1}$ is a sink, there is the 
path $p:u\ra i_{v+1}$ in $Q$. 
As $N$ is an $A$-projective $D$-module and also by Lemma \ref{gradingtool}, 
$$0\not= pN^0_u=(\alpha_{i_{v+1}-1} 
\alpha^*_{i_{v+1}-1})^x N^0_{i_{v+1}}, $$ where $x$ is the cardinality of $[u+1,i_{v+1}]\cap J_v$. 
Similarly,   
$$0\neq pM^{0}_u =(\alpha_{i_{v+1}-1}\alpha^*_{i_{v+1}-1})^y        
M^{0}_{i_{v+1}},$$ 
where $y$ is the cardinality of $[u+1,i_{v+1}]\cap I_v$. Since 
$I_v\subseteq J_v$ we have $x\geq y$. On the other hand, as $f(pN^{0}_u)=pM^{0}_u$, 
we have 
$$(\alpha_{i_{v+1}-1}\alpha^*_{i_{v+1}-1})^y
M^0_{i_{v+1}}\subseteq \mathcal{J}^x M,$$ 
and so $y\geq x$. This shows that $y=x$ and so $[u+1,i_{v+1}]\cap I_v=[u+1,i_{v+1}]\cap J_v$, that is  $(\dimv_\Delta N)_w = 
(\dimv_\Delta M)_w$ for all $u< w \leq i_{v+1}$. If $i_{v+1}=n$
or if $i_{v+1}\not\in J_v$, then we are done. Assume that  $i_{v+1}<n$ 
and $i_{v+1}\in J_v$. Since $J_{v+1}\subseteq I_{v+1}$ and 
$f(N^0_{i_{v+1}})\not \subseteq\mathcal{J}M$ we have 
$I_{v+1}=J_{v+1}$. Then it follows by induction that $I_s=J_s$ for 
all $s\geq v+1$, and so $(\dimv_\Delta N)_w=(\dimv_\Delta M)_w$ 
for all $w>u$.

(b) The proof is similar to (a).
\end{proof}

Let $u=i_v$ be a source and let $M\leq_u N$ be two indecomposable 
rigid $D$-modules supported, and therefore also $\Delta$-supported at $u$.
In this case, $\Hom(M,N)^0_u=\Hom(M,N)_u$, where $\Hom(M,N)_u$
are restrictions of homomorphisms to $M_u$. We have the following lemma, 
similar to Lemma \ref{gradedhom}. Note that $\leq_u$ is defined
on the interval $i_{v-1} < u = i_v$  and so the inequalities are 
opposite to the previous lemma.

\begin{lemma}\label{gradedhomsource}
Let $M$ and $N$ be as above. Then
\begin{itemize}
\item[(a)] $\Hom(N,M)_u=k$ if and only if $(\dimv_\Delta M)_w =
(\dimv_\Delta N)_w$ for all $w<u$.
\item[(b)] $\Hom(M,N)_u=k$ if and only if $(\dimv_\Delta M)_w =
(\dimv_\Delta N)_w$ for all $w>u$.
\end{itemize}
\end{lemma}
%
%
%

\section{Main results}

In the remainder of this paper $Q$ is a quiver obtained by attaching 
a vertex $1$ with an arrow $\gamma$ to an interior vertex $u$ in a quiver 
of type $\mathbb{A}_{n-1}$ with vertices $\{2,\cdots,n\}$ and arrows
between $i$ and $i+1$, and $Q'$ is the quiver with the same underlying graph 
as $Q$ and a unique sink at $u$. Let $P$ be a projective representation of $Q$
with $P=P(d)$ for a vector $d\in \mathbb{N}^ n$. We will study generic 
$\Aut_A P$-orbits in $\Rep(D,P)$, case by case, with respect to the following eight 
orientations of the three arrows incident to vertex $u$.

$$
\xymatrix@=5mm{
& A)  & 1 &     && E)  & 1 \ar[d]^\gamma &     & \\
& u-1 \ar[r]^{\;\;\;\delta} & u \ar[r]^{\beta\;\;\;\;\;} \ar[u]_\gamma & u+1  && u-1 & u \ar[l]_{\;\;\;\;\;\delta} & u+1 \ar[l]_{\beta\;\;\;}& \\
& B)  & 1 &     && F)  & 1 \ar[d]^\gamma &     & \\
& u-1 & u \ar[l]_{\;\;\;\;\;\;\delta} \ar[u]_\gamma \ar[r]^{\beta\;\;\;\;\;} & u+1 && u-1\ar[r]^{\;\;\;\delta} & u & u+1\ar[l]_{\beta\;\;\;} & \\
& C)  & 1\ar[d]^\gamma &   &  & G)  & 1 &     & \\
& u-1 & u \ar[r]^{\beta\;\;\;} \ar[l]_{\;\;\;\;\;\;\delta} & u+1 && u-1 \ar[r]^{\;\;\;\;\delta}& u\ar[u]_\gamma & u+1\ar[l]_{\beta\;\;\;} & \\
& D)  & 1 &     && H)  & 1 \ar[d]^\gamma &     & \\
& u-1 & u \ar[l]_{\;\;\;\;\;\;\delta} \ar[u]_\gamma & u+1 \ar[l]_{\beta\;\;\;}& & u-1 \ar[r]^{\;\;\;\;\delta}& u \ar[r]^{\beta\;\;\;}& u+1 &}
$$

The following theorem is the main technical result of this paper. It shows that 
generic $\Aut_AP$-orbits in $\Rep (D,P)$ can be studied using generic
representations of $Q'$. We emphasize that the theorem is true not only for 
Dynkin quivers, but also for some quivers of tame and even wild type. This 
theorem will be the key step in the proof of the main result stated in the 
introduction. 

\begin{theorem} \label{technicaltheorem}
Let $P$ be a projective representation of $Q$. If the orientation at $u$ is as in
\begin{itemize}
\item[(1)] Case A, B or D, or
\item[(2)] Case C and $u=3,4,n-2$ or $n-1$,
\end{itemize} 
then there is a dimension vector $c$ such that the $\Aut_A P$-space 
$\Rep (D,P)$ is generically equivalent to the $\Gl(c)$-space $\Rep (Q',c)$.
\end{theorem}

The proof of the theorem is long and we postpone the details to next
section. The theorem has the following corollary which follows from
Lemma \ref{openlemma}.

\begin{corollary} \label{technicalcorollary}
Let $Q$, $Q'$ and $P$, $c$ be as in the theorem above. Then there exists a 
dense open $\Aut_A P$-orbit in $\Rep (D,P)$ if and only if there is a 
dense open $\Gl(c)$-orbit in $\Rep (Q',c)$.
\end{corollary}

For a quiver $Q$, let $Q^{op}$ be the opposite quiver of $Q$ with vertices
$Q^{op}_0=Q_0$ and arrows $Q^{op}_1=Q^*_1$, let $A^{op}=kQ^{op}$, 
and note that Case E, F, G and H are opposite to Case A, B, C and D, respectively. 
For each dimension vector $d\in \mathbb{N}^n$, there is a projective 
$A^{op}$-module denoted by $P^{op}=P^{op}(d)$, and an anti-isomorphism 
of algebras 
$$
\phi:\End_{A}P\ra \End_{A^{op}}P^{op}
$$
which maps $\Aut_{A}P$ and $\radE_{A}P$ onto $\Aut_{A^{op}}P^{op}$
and $\radE_{A^{op}}P^{op}$, respectively. 

\begin{lemma} \label{duallemma}
There is an open $\Aut_{A}P$-orbit in $\radE_{A}P$ if and only if there is an 
open $\Aut_{A^{op}}P^{op}$-orbit in $\radE_{A^{op}}P^{op}$.
\end{lemma}
\begin{proof}
The anti-isomorphism $\phi$ induces a commutative diagram 
$$
\xymatrix{\Aut_{A}P\times \radE_{A}P\ar[rr]^{\phi\times \phi}
\ar[d] && \ar[d] \Aut_{A^{op}}P^{op}\times \radE_{A^{op}}P^{op}
\\\radE_{A}P \ar[rr]^{\phi} && \radE_{A^{op}}P^{op}}
$$ 
where the vertical maps are actions and the horizontal maps are 
restrictions of $\phi$. The map $\phi|_{\radE_{A}P}$ sends an open 
$\Aut_{A}P$-orbit in $\radE_A P$ bijectively onto an open 
$\Aut_{A^{op}}P^{op}$-orbit in $\radE_{A^{op}}P^{op}$. 
This proves the lemma.
\end{proof}

The corollary will allow us to reduce from the eight cases A-H to the four
cases A-D to prove Theorem \ref{QuiverTheorem}. We now prove Lemma 
\ref{nonfinite} and Theorem \ref{QuiverTheorem} together, which we 
rewrite as follows.

\begin{theorem}
Let $Q$ be a quiver. Then there is a dense open $\Aut_A P$-orbit in $\radE_A P$ 
for all projective representations $P$ of $Q$ if and only if $Q$ is a Dynkin quiver.
\end{theorem}
\begin{proof}
First assume that $Q$ is a Dynkin quiver and let $P=P(d)$ be a projective
representation of $Q$. By Theorem \ref{JSY} we may assume that $Q$ is 
not of type $\mathbb{A}$ and $Q$ is as in one of the eight cases A-H above.

In cases A-D, using Corollary \ref{technicalcorollary} and Theorem 
\ref{Gabriel}, there is a dense open $\Aut_A P$-orbit in $\Rep (D,P)$, 
and therefore a dense open $\Aut_A P$-orbit in $\radE_A P$, by Theorem 
\ref{equivisom}. Using Lemma \ref{duallemma}, it follows that 
there is a dense open $\Aut_A P$-orbit in $\Rep(D,P)$ in Case E, F, G and H. 

Conversely, assume $Q$ is not a Dynkin quiver. The space $(\radE_A P)^2$ is 
closed under the action on $\Aut_A P$, and so if there is a dense open 
$\Aut_A P$-orbit in $\radE_A P$, then there is a dense open $\Aut_A P$-orbit 
in $\radE_A P/(\radE_A P)^2$. There are a surjective morphism 
$\Aut_AP \rightarrow \Gl(d)$, where the kernel acts trivially on $\radE_A P/(\radE_A P)^2$,
and an isomorphism of vector spaces $\radE_A P/(\radE_A P)^2 \rightarrow \Rep(Q^{op},d)$,
such that there is a commutative diagram
$$\xymatrix{\Aut_A P\times (\radE_A P/(\radE_A P)^2) \ar[d] \ar[r] & \Gl(d) \times \Rep(Q^{op},d) \ar[d] \\  \radE_A P/(\radE_A P)^2 \ar[r] & \Rep(Q^{op},d) }$$
where the vertical maps are actions. 
Since $Q^{op}$ is not Dynkin, there is a dimension vector $d$ such
that there is no dense open $\Gl(d)$-orbit in  $\Rep(Q^{op},d)$. Therefore,
for this $d$, there is no dense open $\Aut_A P(d)$-orbit in $\radE_A P(d)$.
This completes the proof of the theorem.
\end{proof}

\section{Proof of Theorem \ref{technicaltheorem}}

Let $Q$, $Q'$ and $P=P(d)$ for $d\in\mathbb{N}^n$ be as in the 
previous section. We prove Theorem \ref{technicaltheorem} by considering 
Case A-D separately. The proofs are long, but for all four cases we will 
follow the following common strategy. We fix the module structure on the double $\tilde{\Gamma}$
of a subquiver $\Gamma\subseteq Q$ of type $\mathbb{A}$ supported on
$\{1,\cdots,u\}$ in Case C, and on $\{2,\cdots,n\}$ in the other cases. We
then show that the local structure on arrows not in $\tilde{\Gamma}$ is generically 
equivalent to the $\Aut_A P$-space $\Rep(D,P)$, and also generically equivalent 
to the $\Gl(c)$-space $\Rep (Q',c)$ for some dimension vector $c$ constructed 
from $P$. To prove the generic equivalence we use the filtration defined in Section
\ref{grading} and the properties of homomorphisms with respect
to this filtration given in Lemmas \ref{gradedhom} and \ref{gradedhomsource}.

\subsection{Case A}
 
Let $d'$ be the vector given by $(d')_1=0$ and $(d')_i=d_i$ for $i\neq 1$. 
Let $\Gamma$ be the full subquiver of $Q$ supported on $\{2,\cdots,n\}$. 
Since $\Gamma$ is of type $\mathbb{A}$, there is a rigid module 
with $\Delta$-dimension vector $d'$, which we 
denote by $Y(d')$, for the double quiver $\tilde{\Gamma}$ with the corresponding 
relations. We expand $Y(d')$ to an $A$-projective $D$-module $X(d')$ with 
$\Delta$-dimension vector $d'$ as follows. 

On $\tilde{\Gamma}$, $X(d')$ and $Y(d')$ are equal. Let $X(d')_1=Y(d')_u$,  
$X(d')_\gamma=Id$ and $X(d')_{\gamma^*}=X(d')_\delta X(d')_{\delta^*}$, 
where $\delta:u-1\rightarrow u\in Q_1$. Then $X(d')$ is an $A$-projective 
$D$-module, and it is rigid by Theorem \ref{charlemma}, since the construction 
preserves the $\Delta$-dimension vector and the dimension of the endomorphism 
ring.

Let $N^1,\cdots ,N^p$ be indecomposable summands of $Y(d')$, one from each
isomorphism class, supported (but not necessarily $\Delta$-supported) at $u$, 
and ordered such that $N^i<_uN^{i+1}$, where $<_u$ is as in Definition 
\ref{order}. Let $n_i$ denote the multiplicity of $N^i$ as a summand in $Y(d')$,
and let $M^1,\cdots,M^p$ be the corresponding summands of $X(d')$.

As discussed in Section \ref{grading}, $Y(d')$ is graded using the ideal 
$\mathcal{J}$ with the graded components denoted by $Y(d')^i$ and with
$\mathcal{J}Y(d')$ equal to the direct sum of the graded components with positive 
degrees. By construction $X(d')_1=X(d')_u=Y(d')_u=Y(d')^0_u\oplus
\mathcal{J}Y(d')_u$. We denote the two copies of the subspace $Y(d')^0_u$ in 
$X(d')_u$ and $X(d')_1$, by $X(d')^0_u$ and $X(d')^0_1$, respectively.

Following Section \ref{grading}, we let $\End(X(d'))^{0}_u$ denote the 
space of maps which are restrictions 
$$f_{|X(d')_u^0}:X(d')_u^0\rightarrow X(d')_u^0$$ 
of homomorphisms $f:X(d')\rightarrow X(d')$ with 
$f(X(d')^0_u)\subseteq X(d')^0_u$.  The action of $\Aut_D(X(d'))$ on $X(d')$ 
induces by restriction an action of $\Aut(X(d'))^0_u$ on $X(d')^0_u$, 
where $\Aut(X(d'))^0_u \subseteq \End(X(d'))^0_u$ consists of the restrictions
of invertible maps.

\begin{lemma} \label{actionlemma2} The 
$\Aut(X(d'))^0_u\times \Gl_{d_1}$-space $\Hom_k(k^{d_1},X(d')^0_u)$,
with the action given by $(a,g)c=acg^{-1}$ for 
$(a,g)\in \Aut(X(d'))^0_u\times \Gl_{d_1}$ and $c\in \Hom_k(k^{d_1},X(d')^0_u)$,
is a generic section of the $\Aut_A P$-space $\Rep (D,P)$.
\end{lemma}
\begin{proof}
Any $A$-projective $D$-module $X$ with $_{A}X=P$ has a unique submodule 
$X'\subseteq X$ which is generated by the subspaces $X_i$ for $2\leq i \leq n$. 
We have $_{A}X'=P(d')$ and $(X')_{\alpha^*}=X_{\alpha^*}$ for 
$\alpha\in \Gamma$.
We consider the subset in $\Rep (D,P)$ consisting of representations $X$ with 
the submodule $X'=X(d')$. Then as 
$(\tilde{Q}, \mathcal{I})$-representations such $X$ are determined by 
$X(d')\subseteq X$ and the map $X_{\gamma^*}$ as follows, 
$$
\xymatrix{Y(d')^0_u\oplus (\mathcal{J}Y(d'))_u 
\ar@/^1pc/[rrr]^{X_{\gamma}=\tiny\left(\begin{matrix}0 & Id & 0\\ 
0 & 0 & Id \end{matrix}\right)^{tr}} & & & k^{d_1}\oplus Y(d')^0_u\oplus
(\mathcal{J}Y(d'))_u, \ar@/^1pc/[lll]^{X_{\gamma^*}=
\tiny\left(\begin{matrix}c' & 0 & 0 \\ c'' & z_1 & z_2\end{matrix}\right)} } 
$$
where $\left(\begin{matrix}0 & 0 \\ z_1 & z_2\end{matrix}\right) =
X(d')_{\gamma^*}=X(d')_{\delta}X(d')_{\delta^*} $ and $(z_1 \; z_2)$ 
is surjective by the the structure of $Y(d')$. Furthermore, $X$ is determined
by the maps $c'$ and $c''$ and so we denote it by $X(c',c'')$.

Let $$\phi: \Hom_k(k^{d_1},X(d')^0_u)\ra \Rep (D, P)$$ be the morphism given 
by $\phi(c)=X(c, 0)$.  We will show that $\phi$ satisfies the conditions in  
Definition \ref{gensectiondef} of a generic section. First, the subset $\Im \phi$ consists of 
representations of the form $X(c,0)$ and is an affine subspace in $\Rep (D,P)$.
Second, we will show that $\Aut_A P\cdot \Im\phi$ is an open subset 
of $\Rep (D, P)$. 

We first claim that any $X(c', c'')$ is isomorphic to $X(c', 0)$. 
Indeed, a map $$f=(f_i)_{i=1}^n: X(c', c'')\ra X(c, 0),$$ where 
$f_i: X(c', c'')_i\ra X(c, 0)_i$, and $c$ not necessarily equal to $c'$, 
is a homomorphism if and only if the following $3$ conditions are satisfied. 
\begin{itemize}
\item[(1)]$f|_{X(d')}=(a, f_2, \dots, f_n)\in \End_DX(d')$, where 
$a=f_u=f_1|_{X(d')_1}=\left(\begin{tabular}{ll} $a_1$ &$0$ \\ $a_2$& 
$a_3$\end{tabular}\right): X(d')_1\ra X(d')_1$  with respect to the 
decomposition  $X(d')_1=Y(d')^0_u\oplus (\mathcal{J}Y(d'))_u$. 
\item[(2)] $f_1=\left(\begin{tabular}{ll} $g$ &$0$ \\ $b$& $a$
\end{tabular}\right)$, where $g\in \Gl_{d_1}$ and $b: k^{d_1}\ra X(d')_1$.
\item[(3)] (i) $a_1c'=cg$ and  (ii) $a_2c'+a_3c''=(z_1 \; z_2)b$.
\end{itemize}
Here (1) and (2) are due to $X(c', c'')$ and $X(c, 0)$ both having  
$X(d')$ as a submodule, and (3) is because of  
$f_uX(c', c'')_{\gamma^*}=X(c, 0)_{\gamma^*}f_1$. 

When $c=c'$, since $(z_1 \; z_2)$ is surjective, there exists a map 
$b:k^{d_1}\rightarrow X(d')_1$ such that $(z_1 \; z_2)b=c''$. 
If we let $f: X(c',c'')\ra X(c',0)$ be given by $f|_{X(d')}=Id$ and 
$f_1=\left(\begin{matrix}1 & 0 \\ b & 1\end{matrix}
\right)$, then $f$ is an isomorphism.  

Consequently, $\Aut_A P\cdot \Im\phi$ is equal to the subset of 
representations $X$ with the unique submodule $X'\cong X(d')$. Moreover, 
$\Aut_A P\cdot \Im\phi$ is an open subset, because there is a morphism 
$\Rep (D,P)\ra \Rep (D,P(d'))$ which maps $X$ to the unique submodule 
$X'$, with the preimage of the open orbit of $X(d')$ equal to 
$\Aut_A P\cdot \Im\phi$.

Third, if there is an isomorphism $f:X(c',0)\ra X(c,0)$,  then it follows from 
(1)-(3) above that $a_1c'=cg$, i.e. $c, c':  k^{d_1}\ra X(d')^0_u$ are 
conjugate under the action of $\Aut(X(d'))^0_u\times \Gl_{d_1}$ on 
$\Hom_k(k^{d_1}, X(d')^0_u)$. Conversely, suppose that $c$ and $c'$ 
are conjugate using $(a_1,g)\in \Aut(X(d'))^0_u\times \Gl_{d_1}$ and
let $h:X(d')\ra X(d')$ be an automorphism with $a_1$ equal to the restriction
to $X(d')^0_u$. Then since $(z_1 \; z_2)$ is surjective, there exists 
$b:k^{d_1}\rightarrow X(d')_1$ such that $h$, $g$ and $b$ together
give an isomorphism $f:X(c',0)\ra X(c,0)$. This shows that 
$$\phi((\Aut(X(d'))^0_u\times \Gl_{d_1})\cdot c)=
\Aut_A P\cdot \phi(c) \cap \Im \phi,$$ 
for all $c\in \Hom_k(k^{d_1}, X(d')^0_u)$, and so the
 $\Aut(X(d'))^0_u\times \Gl_{d_1}$-space $\Hom_k(k^{d_1},X(d')^0_u)$ 
is a generic section of the $\Aut_A P$-space $\Rep (D,P)$. 
\end{proof}

We compute the maps in $\Aut(X(d'))^0_u$.

\begin{lemma} \label{maplemma} 
Let $i,j\in \{1,\cdots,p\}$. 
\begin{itemize}
\item[(1)] If $\Hom(M^i,M^j)^0_u=k$ then
\begin{itemize}
\item[(a)] $i\geq j$ and $(\dimv_\Delta M^i)_w =
(\dimv_\Delta M^j)_w$ for all $w> u$, or
\item[(b)] $i<j$ and $(\dimv_\Delta M^i)_w =
(\dimv_\Delta M^j )_w$ for all $w\leq u$.
\end{itemize}
\item[(2)] If both (a) and (b) fail then $\Hom(M^i,M^j)^0_u=0$. 
\end{itemize}
\end{lemma}
\begin{proof}
By construction of $M^i$ and $M^j$, we have $\Aut(X(d'))^0_u=\Aut(Y(d'))^0_u$, 
so the lemma follows from Lemma \ref{gradedhom}.
\end{proof}

For each interval 
$[i,j]=\{i,i+1,\cdots,j\}$ for $2\leq i \leq j \leq n$ there is an associated 
indecomposable representation $M[i,j]$ of $Q'$ with support equal to 
$[i,j]$. We construct a representation 
$$Z(d')=\bigoplus_{i=1}^p (Z^i)^{n_i}$$ 
of $Q'$, with $\Aut_{kQ'}(Z(d'))_u\cong \Aut(X(d'))^0_u$, as follows.   
Let $Z^1=M[2,u]$. Given $Z^i=M[j,j']$, let 
$$Z^{i+1}=\left\{\begin{matrix} M[j,j'+1]  \;\;\;\;\;\;& \mbox{ if } \Hom(M^i,M^{i+1})^0_u=k, \\
M[j+1,j']  \;\;\;\;\;\;& \mbox{ if } \Hom(M^{i+1},M^{i})^0_u=k, \\
M[j+1,j'+1] & \mbox{ otherwise.  \;\;\;\;\;\;\;\;\;\; \;\;\;\;\; \;\;\;} \end{matrix} \right.$$

\begin{lemma}
We have $2\leq j\leq u\leq j' \leq n$ for any summand $Z^i=M[j,j']$
in $Z(d')$.
\end{lemma}
\begin{proof}
The inequalities $2\leq j$ and $u\leq j'$ are trivial. 
If $\Hom(M^{i+1},M^{i})^0_u=0$, then by Lemma \ref{maplemma},
$M^i$ and $M^{i+1}$ are $\Delta$-supported differently on
the interval $[u+1,n]$. By the construction of $X(d')$ there are at most
$n-u$ summands $M^i$ with $\Hom(M^{i+1},M^{i})^0_u=0$. Therefore
$j'\leq n$. Similarly, $\Hom(M^{i},M^{i+1})^0_u=0$ for at most
$u-2$ summands $M^i$. Therefore $j\leq u$. 
\end{proof}

Let $c(d')=\dimv Z(d')$ and  
$c(d)=(d_1,c(d')_2,\cdots, c(d')_n).$ 

\begin{lemma} \label{actionlemma}
The $\Aut_{kQ'}(Z(d')_u)\times \Gl_{d_1}$-space $\Hom_k(k^{d_1},Z(d')_u)$ is a
generic section of the $\Gl(c(d))$-space $\Rep (Q',c(d))$.
\end{lemma}
\begin{proof}
Let $\Gamma'$ be the full subquiver of $Q'$ supported on $\{2,\cdots,n\}$.
The $\Aut_{kQ'}(Z(d'))\times \Gl_{d_1}$ orbits in $\Hom_k(k^{d_1},Z(d')_u)$
parametrise the representations of $Q'$ with restriction to
$\Gamma'$ equal to $Z(d')$. The lemma follows if we can prove that
$Z(d')$ is a rigid representation.

The full subquiver $\Gamma'$  has a unique 
sink $u$. Let $M=M[i,j]\oplus M[i',j']$, for $1\leq i,i' \leq u \leq
j,j' \leq n$. Then $\Ext^1_{kQ'}(M[i,j],M[i',j'])\neq 0$ if and only if
$i<u<j$ and $[i',j']\subseteq [i+1,j-1]$. The lemma follows.
\end{proof}

Finally we can prove Case A of Theorem \ref{technicaltheorem}.

\begin{lemma} \label{casea} The 
$\Gl(c(d))$-space $\Rep (Q',c(d))$ is generically equivalent to the 
$\Aut_A P$-space $\Rep (D,P)$.
\end{lemma}
\begin{proof}
There are isomorphisms $(Z^i)_u\rightarrow (M^i)^0_u$, which extend
to an isomorphism $$\Hom_k(k^{d_1},Z(d')_u)\rightarrow
\Hom_k(k^{d_1},X(d')^0_u)$$ of vector spaces. By the construction
$\Hom(M^i,M^j)^0_u\cong \Hom_{kQ'}(Z^i,Z^j)_u$, and we 
have an isomorphism $\Aut_{kQ'}(Z(d'))_u
\rightarrow \Aut(X(d'))^0_u$. Therefore there is a commutative diagram
$$
\xymatrix{Aut_{kQ'}(Z(d'))\times \Gl_{d_1} \times \Hom_k(k^{d_1},Z(d')_u)
\ar[r] \ar[d] &  \Aut(X(d'))^0_u\times \Gl_{d_1} \times
\Hom_k(k^{d_1},X(d')^0_u)\ar[d] \\ \Hom_k(k^{d_1},Z(d')_u) \ar[r] &
\Hom_k(k^{d_1},X(d')^0_u) }
$$
where the vertical maps are actions and the horizontal maps are isomorphisms. 
The lemma now follows from Lemmas \ref{actionlemma2} and \ref{actionlemma}.
\end{proof}

\begin{example} 
Let $Q$ be the quiver
$$\xymatrix@=5mm{ &&&&& 1\\2 \ar[r]&3\ar[r] &4&
 5\ar[l]&6\ar[l]\ar[r]&7\ar[r]\ar[u]&8,}$$ 
$u=7$ and let $d=(2, 2,3,3,0, 3,1,2)$. So $d'=(0,2,3,3,0, 3,1,2)$.
Then $X(d')=M_1\oplus M_2 \oplus M_3$, where $M^1<_7M^2<_7M^3$ 
are as follows. 

$$
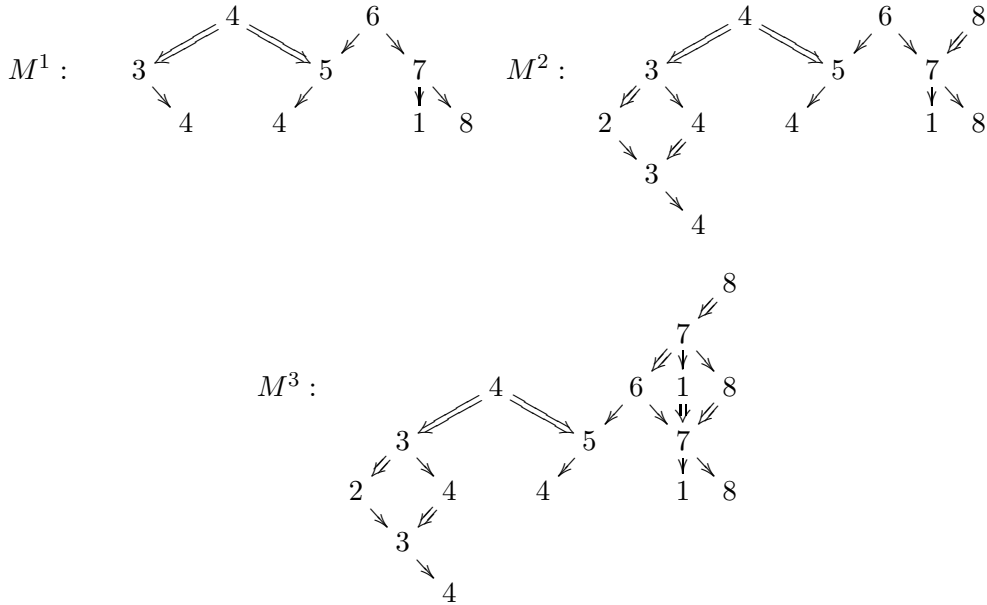
\xymatrix@=2mm{
&&&& 4\ar@{=>}[dll]\ar@{=>}[drr] &&& 6\ar[dr]\ar[dl]       
&&&&&&& 4\ar@{=>}[dll]\ar@{=>}[drr] &&&6\ar[dr]\ar[dl] & & 8\ar@{=>}[dl]
\\ 
M^1: && 3\ar[dr]&&&&5\ar[dl]&&7\ar[d]  \ar[dr]
&&M^2: & &  3\ar[dr]\ar@{=>}[dl]&&&&5\ar[dl] && 7\ar[dr]\ar[d]&&
\\
&&&4&&4&&&1 &8& 
& 2\ar[dr]&&4\ar@{=>}[dl] &&4&&&1 &8&
\\ &&& &&& &&&& & &3\ar[dr] \\
 &&& &&& &&&& && & 4
}$$

$$\xymatrix@=2mm{  
&&&&&&&& & 8\ar@{=>}[dl]\\
&&&&&&&& 7\ar@{=>}[dl]\ar[dr]\ar[d]& \\
M^3:&&&& 4\ar@{=>}[dll]\ar@{=>}[drr] &&& 
6\ar[dl]\ar[dr] &1\ar@{=>}[d]&8\ar@{=>}[dl]
&   \\
&& 3\ar@{=>}[dl]\ar[dr]&&&&5\ar[dl]&&7\ar[d]  \ar[dr] &&\\
&2\ar[dr]&&4\ar@{=>}[dl] &&4&&&1 &8\\
&&3\ar[dr] \\ &&& 4&&& &&
} $$

By computing or by Lemma \ref{maplemma},  
$$\Aut(X(d'))^{0}_7=\{\left(\begin{tabular}{lll} $a_{11}$ & $0$ & $0$ \\ $0$ & $a_{22}$ & $a_{23}$
\\$0$ & $0$ & $a_{33}$\end{tabular}\right) | a_{ij}\in k \}.$$ 

(2) 
We have $$Q': \xymatrix@=5mm{ &&&&& 1 \ar[d]\\ 2 \ar[r]&3\ar[r] &4\ar[r]& 5\ar[r]&6\ar[r]\ar[r]&7&8\ar[l].}$$
We construct the representation $Z(d')=Z^1\oplus Z^2\oplus Z^3$ of $Q'$ with 
$Z^1=M[2,7]$, $Z^2=M[3, 8]$ and $Z^3=[4,8]$. Then 
$\Aut_{kQ'}Z\cong \Aut (X(d'))^{0}_7$. A representation of $Q'$ with dimension
vector $c(d)=(2,1,2,3,3,3,3,2)$ and $Z(d')$ as a submodule is as follows.

$$\xymatrix@=15mm{ &&&&& k^2 \ar[d]^{\tiny c=\left(\begin{matrix}c_{11} & c_{12} 
\\ c_{21} & c_{22} \\ c_{31} & c_{32} \end{matrix}\right)}
\\ k \ar[r]^{\tiny \left(\begin{matrix}1 \\ 0\end{matrix}\right)} & 
k^2\ar[r]^{\tiny \left(\begin{matrix}1 & 0 \\ 0 & 1 \\ 0 & 0\end{matrix}\right)}  &  
k^3\ar[r]^{\tiny \left(\begin{matrix}1 & 0 & 0 \\ 0 & 1 & 0 \\ 0 & 0 & 1\end{matrix}\right)} & 
k^3\ar[r]^{\tiny \left(\begin{matrix}1 & 0 & 0 \\ 0 & 1 & 0 \\ 0 & 0 & 1\end{matrix}\right)} &
k^3\ar[r]^{\tiny \left(\begin{matrix}1 & 0 & 0 \\ 0 & 1 & 0 \\ 0 & 0 & 1\end{matrix}\right)}
&k^3&k^2\ar[l]^{\tiny \left(\begin{matrix}0 & 0 \\ 1 & 0 \\ 0 & 1\end{matrix}\right)}.}$$

A rigid $kQ'$-module with dimension vector $(2,1,2,3,3,3,3,2)$ is isomorphic to
the direct sum of the indecomposable representations with dimension vector 
$(1,1,2,2,2,2,2,1)$, $(1,0,0,1,1,1,1,1)$, respectively, where we may choose
$$c=\left(\begin{matrix}1 & | & 0  \\ 1 & | & 0 \\ - & + & - \\ 0 & | & 1 \end{matrix}\right)$$

(3) The rigid $kQ'$-module in (2) induces a rigid $A$-projective $D$-module
with $\Delta$-dimension vector $d$ as follows. 

$$\xymatrix@=2mm{
&&&&&&& 1 \ar@{=>}[dd]^1\\ &&&&&&&& 8\ar@{=>}[dl]\\
&&&&&&& 7\ar@{=>}[dl]\ar[dr]\ar[d]& \\
&&& 4\ar@{=>}[dll]\ar@{=>}[drr] &&& 6\ar[dl]\ar[dr] &1\ar@{=>}[d]&8\ar@{=>}[dl]\\
& 3\ar@{=>}[dl]\ar[dr]&&&&5\ar[dl]&&7\ar[d]  \ar[dr] & \\
2\ar[dr]&&4\ar@{=>}[dl] &&4&&&1 &8 \\
&3\ar[dr]&&&& &&  \\
&&4 &&& &&}
$$

$$\xymatrix@=2mm{&&&&&&&1\ar@{=>}[dd]^{\tiny \left(\begin{matrix}1 \\ 1\end{matrix}\right)} \\  &&&4\oplus 4 \ar@{=>}[dll]\ar@{=>}[drr]
&&&6\oplus 6\ar[dr]\ar[dl] && 8\oplus 8\ar@{=>}[dl]
\\
&3\oplus 3\ar[dr]\ar@{=>}[dl]_{\tiny \left(\begin{matrix}1 \\ 0\end{matrix}\right)}&&&&5\oplus 5\ar[dl] && 7\oplus 7\ar[dr]\ar[d]&&
\\
2\ar[dr]&&4\oplus 4\ar@{=>}[dl] 
&&4\oplus 4&&&1\oplus 1 &8\oplus 8&\\
& 3\ar[dr]\\
&&4&
}$$
where $j\oplus j$ means that the vector space is $2$-dimensional, 
the matrices on arrows connecting vector spaces of equal dimension are 
identity matrices. Also, the matrices for $1\ra 7$ and $1\ra 7\oplus 7$
are the diagonal blocks of the matrix $c$ in (2).
\end{example}

\subsection{Case B}

This case is very similar to Case A, the main difference being that 
the filtration in Section \ref{grading} is not essential in this section. 
For completeness and the convenience of the reader 
we include most of the details of the proof. 
Let $d'$, $\Gamma$, $Y(d')$, $N^i$ and $n_i$ for $i=1\cdots p$ 
be as in Case A.  We construct a rigid $A$-projective $D$-module $X(d')$ 
with $\Delta$-projective dimension vector $d'$ as follows.
Let $X(d')$ be equal to $Y(d')$ on $\tilde{\Gamma}$. 
Let $X(d')_1=X(d')_u$, $X(d')_\gamma=Id$ and $X(d')_{\gamma^*}=0$ for $\gamma:
u\ra 1\in Q_1$. As in Case A, $X(d')$ is rigid by Theorem \ref{charlemma}. 
Let $M^i$ be the summand in $X(d')$ corresponding to the summand 
$N^i$ in $Y(d')$.

The following lemma is similar to Lemma \ref{actionlemma2} in Case $A$.
Note that since $u$ is a source, we have $X(d')_u=X(d')^0_u$ and 
$\Aut_D(X(d'))_u=\Aut(X(d'))^0_u$, by Lemma \ref{gradingtool}.

\begin{lemma} \label{actionlemmab} The
$\Aut_D(X(d'))_u\times \Gl_{d_1}$-space $\Hom_k(k^{d_1},X(d')_u)$ with the
action $(a,g)c=acg^{-1}$ for $g\in \Gl_{d_1}, a\in \Aut_D(X(d'))_u$ and 
$c\in \Hom_k(k^{d_1},X(d')_u)$, is a
generic section of the $\Aut_A P$-space $\Rep (D,P)$.
\end{lemma}
\begin{proof}
Any $D$-module $X$ with $_{A}X=P(d)$ has a unique submodule $X'\subseteq X$
generated by $X_i$ for all $i\neq 1$. We have $_{A}X'=P(d')$ and
$(X')_{\epsilon^*}=X_{\epsilon^*}$ for $\epsilon\in \Gamma_1$.
We consider the subset of $\Rep (D,P)$ consisting of those representations
$X$ with $X'=X(d')$. So as a $(\tilde{Q},\mathcal{I})$-representation, $X$ is
determined by $X(d')\subseteq X$ and the map $X_{\gamma^*}$ as follows,

$$\xymatrix{Y(d')_u
\ar@/^1pc/[rrr]^{X_{\gamma}=\tiny\left(
\begin{matrix}0 \\ Id \end{matrix}\right)} & & & k^{d_1}\oplus Y(d')_u, \ar@/^1pc/[lll]^{X_{\gamma^*}=\tiny\left(\begin{matrix}c & 0 \end{matrix}\right)} } 
$$ where $c:k^{d_1}\ra Y(d')_u$. 
We let $X(c)$ denote such a representation.

Let $$\phi:\Hom_k(k^{d_1},Y(d')_u)\ra \Rep (D,P)$$ denote the morphism 
$\phi(c)=X(c)$. We show that the $\Aut_D(X(d'))_u\times \Gl_{d_1}$-space \\ 
$\Hom_k(k^{d_1},X(d')_u)$ is a generic section of $\Aut_A P$-space $\Rep (D,P))$, 
using the map $\phi$.

First, $\Aut_A P\cdot \Im\phi$ is equal to the set of representations $X$ with $X'$ 
isomorphic to $X(d')$. There is a morphism $\Rep (D,P)\ra \Rep (D,P(d'))$ 
such that  $\Aut_A P\cdot \Im\phi$ is equal to the preimage of the open orbit
of $X(d')$. Therefore $\Aut_A P\cdot \Im\phi$ is an open subset of 
$\Rep (D,P)$. 

Second, $X(c)\cong X(c')$ if and only if there exists a pair $$(a,\left(\begin{matrix}
g & 0 \\ b & a\end{matrix}\right)),$$ for $a\in \Aut_D(X(d'))_u, g\in \Gl_{d_1},
b:k^{d_1}\ra Y(d')_u$ such that $c'g=ac$, i.e. $c'=acg^{-1}$.
In other words, $X(c)\cong X(c')$ if and only if $c'$ and $c$ are conjugate
under the action of $\Aut_D(X(d'))_u\times \Gl_{d_1}$. 

So the $\Aut_D(X(d'))_u\times \Gl_{d_1}$-space $\Hom_k(k^{d_1},X(d')_u)$ 
is a generic section of the $\Aut_A P$-space $\Rep (D,P)$.
\end{proof}

Similar to Case A, we compute the maps in $\Aut(X(d'))^0_u$. Note that
the inequalities are opposite compared to Case A due to the definition
of $\leq_u$, and $\Hom(M^i,M^j)^0_u=\Hom(M^i,M^j)_u$.

\begin{lemma} \label{map1} 
Let $i,j\in \{1,\cdots,p\}$. 
\begin{itemize}
\item[(1)] If $\Hom(M^i,M^j)^0_u=k$ then
\begin{itemize}
\item[(a)] $i\geq j$ and $(\dimv_\Delta M^i)_w =
(\dimv_\Delta M^j)_w$ for all $w\leq u$, or
\item[(b)] $i<j$ and $(\dimv_\Delta M^i)_w =
(\dimv_\Delta M^j )_w$ for all $w>u$.
\end{itemize}
\item[(2)] If (a) and (b) fail then $\Hom(M^i,M^j)^0_u=0$. 
\end{itemize}
\end{lemma}
\begin{proof}
By construction, we have $\Aut_D(X(d'))_u=\Aut_D(Y(d'))_u$ and so we need 
only prove the corresponding lemma for $N^i$ and $N^j$. Similar to Case A,
the lemma now follows from Lemma \ref{gradedhomsource}.
\end{proof}

Using a similar procedure as in Case $A$, we construct a representation 
$$Z(d')=\oplus (Z^i)^{n_i}$$ of $Q'$, with 
$\Aut_{kQ'}(Z(d'))_u\cong \Aut_D(X(d'))_u$, as follows. 
Let $Z^1=M[u,n]$. Given $Z^i=M[j,j']$, let 
$$Z^{i+1}=\left\{\begin{matrix} M[j-1,j']  \;\;\;\;\;\;& \mbox{ if } \Hom_D(M^i,M^{i+1})_u=k, \\
M[j,j'-1]  \;\;\;\;\;\;& \mbox{ if } \Hom_D(M^{i+1},M^{i})_u=k, \\
M[j-1,j'-1] & \mbox{ otherwise.  \;\;\;\;\;\;\;\;\;\; \;\;\;\;\; \;\;\;\;\;} \end{matrix} \right.$$
Let $c(d')=\dimv Z(d')$ and $c(d)=(d_1,c(d')_2,\cdots, c(d')_n).$ 
Similar to Case $A$ we have the following lemmas.

\begin{lemma}
\begin{itemize}
\item[(1)] For each $Z^i=M[j,j']$, we have $2\leq j\leq u \leq j' \leq n$.
\item[(2)] $Z(d')$ is a rigid representation.
\end{itemize}
\end{lemma}

\begin{lemma} \label{actionlemmaB}
The $\Aut_{kQ'}(Z(d')_u)\times \Gl_{d_1}$-space $\Hom_k(k^{d_1},Z(d')_u)$ is a
generic section of the $\Gl(c(d))$-space $\Rep (Q',c(d))$.
\end{lemma}

As in Case A, using Lemma \ref{actionlemmab}, Lemma \ref{map1} and 
Lemma \ref{actionlemmaB}, we have the following lemma which proves Case B 
of Theorem \ref{technicaltheorem}.

\begin{lemma} \label{caseb} 
The $\Gl(c(d))$-space $\Rep (Q',c(d))$ is generically equivalent to the $\Aut_A P$-space
$\Rep (D,P)$.
\end{lemma}

\subsection{Case C}

Unlike case $A$ and $B$, in this case we let $\Gamma$ be the full subquiver of 
$Q$ supported on the vertices $\{1,\cdots,u\}$. 
Moreover, we assume for now that $u=n-2$, and postpone 
the cases where $u=3,4$ or $n-1$ to the proof of Lemma \ref{casec} below.
Let $\alpha$ and $\beta$ be the arrows in $Q$ as follows 
$$\xymatrix@=5mm{& & 1\ar[d]^\gamma & & \\ \cdots & u-1 & 
u\ar[l] \ar[r]^\beta & u+1 \ar[r]^\alpha & \ar[l] u+2}$$
where $\alpha$ could be of either orientation, that is, $u+2$ is either a sink 
or a source.  Let $d'$ be given by $d'_i=0$ for $i=u+1,u+2$ and $d'_i=d_i$ otherwise. 
Let $Y(d')$ be a rigid module for the double quiver $\tilde{\Gamma}$ with the 
corresponding relations, which we extend to a rigid $A$-projective $D$-module 
$X(d')$ as follows. 

On $\tilde{\Gamma}$ we let $X(d')$ be equal to $Y(d')$. Let $X(d')_{u+1}=Y(d')_u$,
$X(d')_\beta=Id$, and $X(d')_{\beta^*}=Y(d')_{\gamma}Y(d')_{\gamma^*}$. 
If $u+2$ is a sink, then $X(d')_{u+2}=Y(d')_u$, $X(d')_\alpha=Id$ and
$X(d')_{\alpha^*}=X(d')_{\beta^*}$, and if $u+2$ is a source, then $X(d')_{u+2}=0$. 
The extension of $Y(d')$ to $X(d')$ preserves the $\Delta$-dimension vector and 
the dimension of the endomorphism ring, and so $X(d')$ is rigid by Theorem 
\ref{charlemma}.

Let $N^1,\cdots ,N^p$ be indecomposable summands of $Y(d')$, one from each 
isomorphism class, supported at $u$ and ordered such that  $N^i<_uN^{i+1}$, and 
let $n_i$ be the multiplicity of $N^i$ as a summand in $Y(d')$. For each 
indecomposable summand $N^i$ of $Y(d')$, let $M^i$ be the corresponding
indecomposable summand of $X(d')$.  

Let $d''$ be the dimension vector supported on $\{u+1,u+2\}$, given by $d''_i=0$ 
for $i\leq u$ and $d''_i=d_i$ otherwise. Let $X(d'')$ be the rigid $D$-module with 
$\Delta$-dimension vector $d''$.

We will construct two groups $H_V$ and $H_W$, which act on vector spaces 
$V$ and $W$ in such a way that the $H_V\times H_W$-space $\Hom_k(V,W)$ is a 
generic section of the $\Aut_A P$-space $\Rep (D,P)$. The construction depends  
on whether $u+2$ is a sink or a source.

If $u+2$ is a sink, let $$V=soc(X(d''))\subseteq X(d'')_{u+2}, \; W=X(d')^0_{u+1},$$ and 
$$H_V=\{f|_V \mid f\in \Aut_D(X(d'')) \}, \; H_W=\Aut(X(d'))^0_{u+1}.$$ 
If $u+2$ is a source, let
$$V=(X(d'')/\rad X(d''))_{u+1}=k^{d_{u+1}}, \; W=X(d')^0_{u},$$ and 
$$H_V = \{ \overline{f}:(X(d'')/\rad X(d''))_{u+1}\rightarrow
(X(d'')/\rad X(d''))_{u+1}\mid f\in \Aut_D(X(d''))\},$$ $$H_W=\Aut(X(d'))^0_{u}.$$
In either case, we have the following lemma.

\begin{lemma} \label{l1}
The $H_V\times H_W$-space $\Hom_k(V,W)$ is generically equivalent to the 
$\Aut_A P$-space $\Rep (D,P)$.
\end{lemma}
\begin{proof}
Associated to any $A$-projective $D$-module $X$ there is the unique submodule 
$X'\subseteq X$ generated by the spaces $X_i$ for $i\leq u$, with corresponding 
quotient denoted by $X''=X/X'$. If $d=\dimv_\Delta X$, then $\dimv_\Delta X'=d'$ 
and $\dimv_\Delta X''=d''$.  We consider the subset of $\Rep (D,P)$ consisting 
of $D$-modules $X$ with $X'=X(d')$ and $X''=X(d'')$, which is nonempty since
it contains the direct sum $X(d')\oplus X(d'')$. Moreover, we fix a decomposition 
of vector spaces $X_i=X(d')_i\oplus X(d'')_i$ for all vertices $i$.

We decompose $$X(d')=M\oplus N\oplus L,$$ where the indecomposable summands of 
$M$ are $\Delta$-supported at both $1$ and $u$, those of  
$N$ are $\Delta$-supported at either $1$ or $u$, but not both, and $L$ is not supported 
at $u$. Also, we decompose 
$$
X(d'')=R\oplus T,
$$ 
where the indecomposable summands of $R$ are $\Delta$-supported at both $u+1$ and 
$u+2$, and those of $T$ are not $\Delta$-supported at both $u+1$ and
$u+2$.

The proof is divided into two parts, depending on whether $u+2$ is a source
or a sink. \\

{\bf Part 1}: We first consider the case where $u+2$ is a sink. \\

Let $W_1=(M^0)_{u+1}$, $W_2=(M^1)_{u+1}$ and $W_3=(N^0)_{u+1}$,
where $(N^0)_{u+1}=N_{u+1}$
We have 
$$
W=W_1\oplus W_3.
$$  
Let $V_1=(R^0)_{u+2}$, $V_2=(R^1)_{u+2}$ and $V_3=T_{u+2}$. 
We have $$V=V_2\oplus V_3$$

By the relations of $D$, such a $(\tilde{Q},\mathcal{I})$-representation
$X$ is determined by $X(d')$, $X(d'')$ and the maps between the vertices $u+1$ 
and $u+2$, which have the form 
$$
\xymatrix@=15mm{X(d')_{u+1} \oplus X(d'')_{u+1} \ar@/^1pc/[rr]^{X_{\alpha}=
\tiny\left(\begin{matrix}Id & 0 \\ 0 & X(d'')_\alpha\end{matrix}\right)} & & 
X(d')_{u+2} \oplus X(d'')_{u+2} \ar@/^1pc/[ll]^{X_{\alpha^*}=
\tiny\left(\begin{matrix}X(d')_{\alpha^*} & c \\ 0 & X(d'')_{\alpha^*} 
\end{matrix}\right)}}$$ where $X(d')_{u+1}=X(d')_{u+2}=
W_1\oplus W_2 \oplus W_3$, $X(d'')_{u+2}=V_1\oplus V_2 \oplus V_3$,
$X(d'')_{u+1}=V_2\oplus V_3$ if $d_{u+1}>d_{u+2}$, and
$X(d'')_{u+1}=V_2$ if $d_{u+1}\leq d_{u+2}$, and $$c=(c_{ij})_{ij}:V_1 
\oplus V_2 \oplus V_3\ra W_1\oplus W_2 \oplus W_3.$$ 
We let $X(c)$ denote the representation in the diagram above.

Our first step is to show that with some restriction on the map $c_{12}$,
$X(c)$ is isomorphic to $X(c_0)$, where  
$$c_0=\left(\begin{matrix}0 & c_{12} & c_{13} \\ 0 & 0 & 0 \\ 0 & c_{32} 
& c_{33} \end{matrix}\right).$$

By fixing a basis we may assume we have matrices 
$$X(d'')_\alpha=\left(\begin{matrix}
0 & 0 \\ Id & 0 \\ 0 & Id\end{matrix}\right), 
X(d'')_{\alpha^*}=\left(\begin{matrix}Id & 0 & 0 \\ 0 & 
0 & 0\end{matrix}\right)  \mbox{ and }  
X(d')_{\alpha^*}=\left(\begin{matrix}0 & 0 & 0 \\ Id & 0 & 0 \\ 0 & 0 &
0\end{matrix}\right)$$ Note that if
$d_{u+1}\leq d_{u+2}$ the second column of $X(d'')_\alpha$ and the
second row of $X(d'')_{\alpha^*}$ are empty.

We choose a basis of $W_1$ (and of $W_2$) such that 
there is one basis element from each summand $M^i$, and order the basis elements according
to the order $\leq_u$.

A map $f=(f_i)_{i=1}^n: X(c)\ra X(c')$, where $f_i: X(c)_i\ra X(c')_i$, 
is a homomorphism if and only if the following $4$ conditions are satisfied. 
\begin{itemize}
\item[(1)]$f|_{X(d')}=(f_1, f_2, \dots, f_u,a,a)\in \End_DX(d')$.
\item[(2)] $(f_{u+1},f_{u+2})=
\left(\left(\begin{matrix}a & b' \\ 0 & g'\end{matrix}\right),
\left(\begin{matrix}a & b \\ 0 & g\end{matrix}\right)\right)$, where 
$b=(b_{ij})_{ij}:V_1\oplus V_2\oplus V_3\ra W_1\oplus W_2\oplus W_3$,
and $b'=b|_{X(d'')_{u+1}}$.
\item[(3)] $\overline{f}=(0, 0,\dots,0,g',g)\in \End_DX(d'')$.
\item[(4)] $c'g=ac-X(d')_{\alpha^*}b+b'X(d'')_{\alpha^*}$, as $f_{u+1}X(c)_{\alpha^*}
=X(c')_{\alpha^*}f_{u+2}$
\end{itemize}

Moreover, as matrices with respect to the fixed bases, we have 
$$-X(d')_{\alpha^*}b+b'X(d'')_{\alpha^*}=
\left(\begin{matrix}b_{12} & 0 & 0 \\ -b_{11}+b_{22} & -b_{12} & -b_{13} 
\\ b_{32} & 0 & 0 \end{matrix}\right) \mbox{  for } b=(b_{ij})_{ij}, $$
$$g=\left(\begin{matrix}g_1 & 0 & 0 \\ x_1 & g_1 & x_2 \\ x_3 & 0 & g_2
\end{matrix}\right)$$ for any invertible matrices $g_1$ and $g_2$ 
and any matrices $x_1,x_2$ and $x_3$, and $$a=\left(\begin{matrix}
a_1 & 0 & 0 \\  y_1 & a_1 & y_2 \\ y_3 & 0 & a_2\end{matrix}\right),$$
where $a_1$ and $a_2$ are invertible. In general, not all matrices 
occur as matrices $a_i$ and $y_i$ of an automorphism $a$.

The groups $H_V$ and $H_W$ consist of matrices 
of the form $$\left(\begin{matrix} 
g_1 & x_2 \\ 0 & g_2\end{matrix}\right) \mbox{ and } \left(\begin{matrix}
a_1 & 0 \\ y_3 & a_2\end{matrix}\right),$$ respectively.

For a given map $c$ we let $b(c)$ be $$b(c)=\left(\begin{matrix}c_{21} 
& c_{22} & c_{23} \\ 0 & 0 & 0 \\ 0 & -c_{31} & 0 \end{matrix}\right).$$ 
Then for any $X(c)$, using an isomorphism $f$ with $a=Id$, $g=Id$ and 
$b=b(c)$, we see that $X(c)\cong X(c')$ for $$c'=\left(\begin{matrix}c'_{11} 
& c_{12} & c_{13} \\ 0 & 0 & 0 \\ 0 & c_{32} & c_{33} \end{matrix}\right),$$ 
where $c'_{11}=c_{11}+c_{22}$. 

Using an isomorphism $f$ with $b=0$, $$a=\left(\begin{matrix}Id & 0 & 0 \\ y & Id & 0 \\ 0 & 0 & Id
\end{matrix}\right):W_1\oplus W_2\oplus W_3\ra W_1\oplus W_2
\oplus W_3,$$ and $$g=\left(\begin{matrix}Id & 0 & 0 \\ x & Id & 
0 \\ 0 & 0 & Id\end{matrix}\right):V_1\oplus V_2\oplus V_3\ra V_1\oplus V_2\oplus V_3,$$  
we have $X(c')\cong X(c'')$, where 
$$c''=ac'g^{-1}=\left(\begin{matrix}c'_{11}-c_{12}x & c_{12} & c_{13} \\ 
y(c'_{11}-c_{12}x) 
& yc_{12} & yc_{13} \\ -c_{32}x & c_{32} & c_{33}\end{matrix}\right).$$ Now by 
an isomorphism with $a=Id$, $g=Id$ and an $b=b(c'')$ as above, $X(c'')$ is isomorphic 
to $X(c''')$ with 
$$c'''=\left(\begin{matrix}c'_{11}-c_{12}x + yc_{12}& c_{12} & c_{13} \\ 
0& 0 & 0 \\ 0 & c_{32} & c_{33}\end{matrix}\right).$$

We identify $W_1$ with $W_2$ and $V_1$ with $V_2$. 
We claim that with some restriction on $c_{12}$, there exists maps $x:W_1\ra W_2$ 
and $y:V_1\ra V_2$ such that $c'_{11}-c_{12}x + yc_{12}=0$. Indeed, by Example 
\ref{gensecexample}, the matrix space $\Mat_{\dim_k W_1\times \dim_k V_1}$ has an open 
$\Gl_{\dim_k V_1}\times \B_{\dim_k W_1}$-orbit,  
where $\B_{\dim_k W_1}$ consists of all invertible
lower triangular $\dim_k W_1\times \dim_k W_1$-matrices, and the action is given by 
$(h,b)c_{12}=bc_{12}h^{-1}$. For a matrix $c_{12}$ with an open 
$\B_{\dim_k W_1}\times \Gl_{\dim_k V_1}$-orbit, the associated map on tangent 
spaces $$(*) \; \; \; (\overline{b},\overline{h})\mapsto \overline{b}c_{12}-c_{12}\overline{h}$$ 
is surjective.  Now by the structure of $X(d'')$, with respect to the fixed basis, the map $x:V_1\ra V_2$ 
can be any quadratic matrix. Also, by the structure of $X(d')$ and by the choice of basis of $W_1$
and $W_2$, any lower triangular quadratic matrix
occur as a matrix of $y:W_1\ra W_2$. So there exists $c_{12}$ such that the
map $(*)$ is surjective. Thus there exist $x$ and $y$ such that 
$c'_{11}-c_{12}x + yc_{12}=0$, which proves the claim. 

Consequently, $X(c)$ is isomorphic to 
$X(c_0)$ where $$c_0=\left(\begin{matrix}0 & c_{12} & c_{13} \\ 0 & 0 & 0 \\ 
0 & c_{32} & c_{33}\end{matrix}\right).$$ 

Let $\iota:\Hom_k(V,W)\ra \Hom_k(V_1\oplus V_2\oplus V_3,W_1\oplus W_2\oplus W_3)$ 
be the map $$\iota(\left(\begin{matrix} z_{12} & z_{13} \\ z_{32} & z_{33} 
\end{matrix}\right))=\left(\begin{matrix} 0 & z_{12} & z_{13} \\ 0 & 0 & 0 \\
0 & z_{32} & z_{33} \end{matrix}\right).$$ 
Let $$\phi:\Hom_k(V, W)\ra \Rep (D,P)$$ be given by $\phi(z)=X(\iota(z))$. 
We first show that $\Aut_A P\cdot \Im\phi$ contains a non-empty subset of 
$\Rep (D,P)$. Let $F$ be the set of all representations of the form $X(c)$
for $c\in \Hom_k(V_1\oplus V_2\oplus V_3,W_1\oplus W_2\oplus W_3)$.
The set $\Aut_A P\cdot F$ is equal to the subset of $\Rep (D,P)$ consisting
of $X$ with $X'\cong X(d')$ and $X''\cong X(d'')$. Using the morphism
of varieties $\Rep (D,P)\ra \Rep (D,P(d'))\times \Rep (D,P(d''))$ 
given by $X\mapsto (X',X'')$, we have that $\Aut_A P\cdot F$ is equal
to the preimage of $\Aut_A P(d')\cdot X(d')\times \Aut_A P(d'')
\cdot X(d'')$ which is open. 
The condition that the map $(*)$ is surjective for $c_{12}$ is an open condition on 
the set of representations in $F$. Therefore $(\Aut_A P\cdot \Im\phi)\cap F$
must contain an open subset of $F$, and so finally $\Aut_A P\cdot \Im\phi$
contains an open subset of $\Rep (D,P)$.

Assume that $X(\iota(z))\cong X(\iota(z'))$. Then there exists an isomorphism 
$f$ with maps 
$$a=\left(\begin{matrix} a_1 & 0 & 0 \\  y_1 & a_1 & y_2 \\ y_3 & 0 & 
a_2\end{matrix} \right) \mbox{ and } g=\left(\begin{matrix}g_1 & 0 & 0 
\\ x_1 & g_1 & x_2 \\ x_3 & 0 & g_2 \end{matrix}\right)$$ 
and some $b$ such that conditions (1)-(4) are satisfied. Then by explicitly computing the
matrices in (4), we see that $z$ and 
$z'$ are conjugate under the action of $H_V\times H_W$ using $$\left(\left(\begin{matrix} g_1 & x_2 \\ 0 & g_2\end{matrix}
\right), \left(\begin{matrix} a_1 & 0 \\ y_3 & a_2\end{matrix}\right)\right).$$ 
Conversely, if $z$ and $z'$ are conjugate under the action of $H_V\times H_W$ 
via
$$\left(\left(\begin{matrix} g_1 & x_2 \\ 0 & g_2\end{matrix}\right), \left(
\begin{matrix} a_1 & 0 \\ y_3 & a_2\end{matrix}\right)\right),$$ 
then $X(\iota(z))\cong X(\iota(z'))$ using an isomorphism $f$ with $b=0$, 
$$a=\left(\begin{matrix} a_1 & 0 & 0 \\  0 & a_1 & 0 \\ 
y_3 & 0 & a_2\end{matrix} \right) \mbox{ and } 
g=\left(\begin{matrix}g_1 & 0 & 0 \\ 
0 & g_1 & x_2 \\ 0 & 0 & g_2 \end{matrix}\right).$$ 

Therefore, the $H_V\times H_W$-space $\Hom_k(V,W)$ is generically equivalent
to the 
the $\Aut_A P$-space in $\Rep (D,P)$, in the case where $u+2$ is a sink.\\

{\bf Part 2:}  We now consider the case where $u+2$ is a source. \\

Let $W_1=(M^0)_{u}$, $W_2=(M^1)_{u}$, $W_3=(N^0)_{u}$ and
$W=W_1\oplus W_3.$  Let $V_1=(R^0)_{u+1}$, $V_2=(R^1)_{u+1}$ 
and $V_3=T_{u+1}$. We have 
$V=V_1\oplus V_3$ if $d_{u+1}>d_{u+2}$ and $V=V_1$ if 
$d_{u+1}\leq d_{u+2}$. 

Recall that the representation $X$ has a submodule $X(d')$ with quotient isomorphic
to $X(d'')$, and that we have a decomposition of vector spaces $X_i = X(d')_i \oplus X(d'')_i$.
By the relations of $D$, we see that as a $(\tilde{Q},\mathcal{I})$-representation
$X$ is determined by $X(d')$, $X(d'')$ and the maps between the vertices 
$u$ and $u+1$ which have the form $$\xymatrix{X(d')_{u} 
\ar@/^1pc/[rr]^{X_{\beta}=\tiny\left(\begin{matrix}Id \\ 0 
\end{matrix}\right)} & & X(d')_{u+1} \oplus X(d'')_{u+1}
\ar@/^1pc/[ll]^{X_{\beta^*}=\tiny\left(\begin{matrix}X(d')_{\beta^*} & c 
\end{matrix}\right)}}$$ where $X(d')_{u}=X(d')_{u+1}=
W_1\oplus W_2 \oplus W_3$, $X(d'')_{u+1}=V_1\oplus V_2 \oplus V_3$, 
and $$c=(c_{ij})_{ij}:V_1 \oplus V_2 \oplus V_3\ra W_1\oplus W_2 \oplus W_3.$$ 
By fixing a basis we may assume $$X(d')_{\beta^*}=
\left(\begin{matrix}0 & 0 & 0 \\ Id & 0 & 0 \\ 0 & 0 & 0\end{matrix}\right).$$ 
Due to the relation $\beta^*\alpha=0$ in $D$ we have $c|_{V_2}=0$ 
and so $$c=\left(\begin{matrix}c_{11} & 0 & c_{13} \\ c_{21} & 0 & c_{23} \\ c_{31} 
& 0 & c_{33}\end{matrix}\right).$$ Note that if $d_{u+1}\leq  d_{u+2}$
then $V_3=0$ and so $c_{i3}=0$ for all $i$.  
Denote by $X(c)$ the representation of the form above determined by
$c=(c_{ij})_{ij}$, where $c_{i2}=0$ for all $i$.

A map $f=(f_i)_{i=1}^n: X(c)\ra X(c')$, where $f_i: X(c)_i\ra X(c')_i$, 
is a homomorphism if and only if the following four conditions are satisfied. 
\begin{itemize}
\item[(i)] $f|_{X(d')}=(f_1, f_2, \dots, f_u,a,0)\in \End_DX(d')$.
\item[(ii)] $(f_{u},f_{u+1})=\left(a,
\left(\begin{matrix}a & b \\ 0 & g\end{matrix}\right)\right)$, where 
$b=(b_{ij})_{ij}:V_1\oplus V_2\oplus V_3\ra W_1\oplus W_2\oplus W_3$.
\item[(iii)] $\overline{f}=(0,0,\dots,0,g,g')\in \End_D(X(d''))$.
\item[(iv)] $c'g=ac-X(d')_{\beta^*}b$ as $f_uX(c)_{\beta^*}= X(c')_{\beta^*}f_{u+1}$.
\end{itemize}

Moreover, as matrices we have $$-X(d')_{\beta^*}b=
\left(\begin{matrix}0 & 0 & 0 \\ -b_{11} & -b_{12} & -b_{13} 
\\ 0 & 0 & 0 \end{matrix}\right) \mbox{ for } b=(b_{ij})_{ij},$$  
$$g=\left(\begin{matrix}g_1 & 0 & 0 \\ x_1 & g_1 & x_2 \\ x_3 & 0 & g_2
\end{matrix}\right),$$ where $g_1$ and $g_2$ can be any invertible
matrices and $x_1$, $x_2$ and $x_3$ can be any matrices, and
$$a=\left(\begin{matrix}
a_1 & 0 & 0 \\  y_1 & a_1 & y_2 \\ y_3 & 0 & a_2\end{matrix}\right).$$
As in Part 1, the matrices $a_i$ and $y_i$ depend on the structure of $X(d')$,
and so not all matrices occur.  
Now $H_V$ is the group of matrices $$\left(\begin{matrix}
g_1 & 0 \\ x_3 & g_2\end{matrix}\right)$$ if $d_{u+1}>d_{u+2}$,
and $g_1$, otherwise. The group $H_W$ consists of matrices 
$$ \left(\begin{matrix}
a_1 & 0 \\ y_3 & a_2\end{matrix}\right),$$ induced by automorphisms
of $X(d')$.

For a given map $c=(c_{ij})_{ij}$ with $c_{i2}=0$ for any $i$, we let $b(c)$ be $$b(c)=\left(\begin{matrix}c_{21} 
& 0 & c_{23} \\ 0 & 0 & 0 \\ 0 & 0 & 0 \end{matrix}\right).$$ 
Then for any $X(c)$, using an isomorphism $f$ with $a=Id$, $g=Id$ and 
$b=b(c)$, we see that $X(c)\cong X(c')$ for $$c'=\left(\begin{matrix}c_{11} 
& 0 & c_{13} \\ 0 & 0 & 0 \\ c_{31} & 0 & c_{33} \end{matrix}\right).$$

Let $\iota:\Hom_k(V,W)\ra \Hom_k(V_1\oplus V_2\oplus V_3,W_1\oplus W_2\oplus W_3)$ 
be the map $$\iota(\left(\begin{matrix} z_{11} & z_{13} \\ z_{31} & z_{33} 
\end{matrix}\right))=\left(\begin{matrix} z_{11} & 0 & z_{13} \\ 0 & 0 & 0 \\
z_{31} & 0 & z_{33} \end{matrix}\right)$$ if $d_{u+1}>d_{u+2}$ and 
$$\iota(\left(\begin{matrix} z_{11} \\ z_{31} 
\end{matrix}\right))=\left(\begin{matrix} z_{11} & 0 & 0 \\ 0 & 0 & 0 \\
z_{31} & 0 & 0 \end{matrix}\right)$$ otherwise.
Let $$\psi:\Hom_k(V, W)\ra \Rep (D,P)$$ be given by $\psi(z)=X(\iota(z))$. 
First, $\Aut_A P\cdot \Im\phi$ contains the nonempty subset of $\Rep (D,P)$ 
consisting of representations $X$ such that $X'\cong X(d')$ and $X''\cong X(d'')$,
and so similar to the case where $u+2$ is a sink, 
$\Aut_A P\cdot \Im\phi$ contains a nonempty  open subset of $\Rep (D,P)$.

Assume that $X(\iota(z))\cong X(\iota(z'))$. Then there exists an isomorphism 
$f$ with maps $b$, 
$$
a=\left(\begin{matrix} a_1 & 0 & 0 \\  y_1 & a_1 & y_2 \\ y_3 & 0 & 
a_2\end{matrix} \right) \mbox{ and } g=\left(\begin{matrix}g_1 & 0 & 0 \\ 
x_1 & g_1 & x_2 \\ x_3 & 0 & g_2 \end{matrix}\right)
$$ 
such that conditions (i)-(iv) are satisfied.  Then by explicitly computing the matrices in (iv), 
we see that $z$ and 
$z'$ are conjugate under the action of $H_V\times H_W$ via
$$
\left(\left(\begin{matrix} g_1 & 0 \\ x_3 & g_2\end{matrix}
\right), \left(\begin{matrix} a_1 & 0 \\ y_2 & a_2\end{matrix}\right)\right).
$$ 
Conversely, if $z$ and $z'$ are conjugate under the action of $H_V\times H_W$ 
using 
$$
\left(\left(\begin{matrix} g_1 & 0 \\ x_3 & g_2\end{matrix}\right), \left(
\begin{matrix} a_1 & 0 \\ y_2 & a_2\end{matrix}\right)\right),
$$ 
then $X(\iota(z))\cong X(\iota(z'))$ using a map $f$ with $b=0$, 
$$
a=\left(\begin{matrix} a_1 & 0 & 0 \\  0 & a_1 & 0 \\ 
y_2 & 0 & a_2\end{matrix} \right) \mbox{ and } g=\left(\begin{matrix}g_1 
& 0 & 0 \\ 0 & g_1 & 0 \\ x_3 & 0 & g_2 \end{matrix}\right).
$$ 

Therefore the $H_V\times H_W$-space $\Hom_k(V,W)$ is generically equivalent
to the $\Aut_A P$-space $\Rep (D,P)$.
\end{proof}

We compute $\Aut(X(d'))^0_u$.

\begin{lemma} \label{maplemmaC} 
Let $i,j\in \{1,\cdots,p\}$. 
\begin{itemize}
\item[(1)] If $\Hom(M^i,M^j)^0_u=k$ then
\begin{itemize}
\item[(a)] $i\geq j$ and $(\dimv_\Delta M^i)_w =
(\dimv_\Delta M^j)_w$ for all $w<u$, or
\item[(b)] $i<j$ and $(\dimv_\Delta M^i)_w =
(\dimv_\Delta M^j )_w$ for all $w=1, u$.
\end{itemize}
\item[(2)] If both (a) and (b) fail then $\Hom(M^i,M^j)^0_u=0$. 
\end{itemize}
\end{lemma}
\begin{proof}
By relabeling the vertices $1\mapsto 1$, $u\mapsto 2$ and $u-t\mapsto t+2$ 
we are in Case A, and the lemma follows from Lemma \ref{maplemma}.
\end{proof}

Similar to Case A, we construct a rigid 
representation $Z(d')$ of $Q'$ with summands supported on intervals in 
$$\xymatrix@=5mm{2  \ar[r] & 3 \ar[r] & \cdots \ar[r] & u-1 \ar[r] & u & \ar[l] 1}.$$ 
We have $$Z(d')=\oplus (Z^i)^{n_i}$$ as follows, where $n_i$ 
is the multiplicity of $M^i$ in $X(d')$.  Let $Z^1=M[u,1]$, that is $Z^1$
is supported on $1$ and $u$. Given $Z^i=M[j,j']$, let
$$Z^{i+1}=\left\{\begin{matrix}\; M[j-1,j'] & \mbox{ if } \Hom(M^i,M^{i+1})^0_u=k, \\
M[j,u]  \;\;\;\;\;\;& \mbox{ if }\Hom(M^{i+1},M^{i})^0_u =k, \\
M[j-1,u] & \mbox{ otherwise.  \;\;\;\;\;\;\;\;\;\; \;\;\;\;\; \;\;\;} \end{matrix} \right.$$
Only one of the two latter cases can occur,
and it occurs at most once, and so every $Z^i$ is supported at $u$.
Let $$c(d')=\dimv Z(d').$$ 

Recall $X(d'')=R\oplus T$, where the indecomposable summands of $R$ has 
$\Delta$-support at both $u+1$ and $u+2$, and $T$ has $\Delta$-support
at $u+1$ or $u+2$, but not both. Let $a$ and $b$ be the multiplicities of the 
indecomposable summand in $R$ and $T$, respectively.
Similarly, we have a rigid representation $Z(d'')$ of $Q'$ as follows. 

If $u+2$ is a 
sink, let $$Z(d'')=M[u+1,u+2]^a\oplus M[u+1,u+1]^b.$$ If $u+2$ is a source,
let $$Z(d'')=M[u+1,u+1]^a\oplus M[u+1,u+2]^b \mbox{ when } d_{u+1}>d_{u+2},$$ 
and $$Z(d'')=M[u+1,u+1]^a \mbox{ when } d_{u+1}\leq d_{u+2}.$$ Let 
$$c(d'')=\dimv Z(d'') \mbox{ and } c(d)=c(d')+c(d'').$$ 

\begin{lemma} \label{l2}
The $\Aut_{kQ'}(Z(d''))\times \Aut_{kQ'}(Z(d'))$-space 
$\Hom_k(Z(d'')_{u+1},Z(d')_{u})$ is a 
generic section of the $\Gl(c(d))$-space $\Rep (Q',c(d))$.
\end{lemma}
\begin{proof}
The $\Aut_{kQ'}(Z(d''))\times \Aut_{kQ'}(Z(d'))$-orbits in $\Hom_k(Z(d'')_{u+1},Z(d')_{u})$
parameterise the representations in $\Rep (Q',c(d))$ with restriction
to the subquiver on $\{1,\cdots u\}$ equal to $Z(d')$ 
and restriction to the subquiver on $\{u+1,
u+2\}$ equal to $Z(d'')$. Clearly $Z(d'')$ is rigid and by the proof 
of Lemma \ref{actionlemma}, $Z(d')$ is rigid. So
the $\Aut_{kQ'}(Z(d''))\times \Aut_{kQ'}(Z(d'))$-space 
$\Hom_k(Z(d'')_{u+1},Z(d')_{u})$ is a 
generic section of the $\Gl(c(d))$-space $\Rep (Q',c(d))$.
\end{proof}

We can now prove Theorem \ref{technicaltheorem} in case C).

\begin{lemma} \label{casec}
Let $P$ be a projective representation of $Q$. If the orientation 
at $u$ is as in Case C and $u=3,4,n-2$ or $n-1$,
then there is a dimension vector $c$ such that the 
$\Aut_A P$-space $\Rep (D,P)$ is generically 
equivalent to the $\Gl(c)$-space $\Rep (Q',c)$.
\end{lemma}
\begin{proof}
In the case $u=n-2$ the proof follows from 
Lemma \ref{l1}, \ref{maplemmaC} and \ref{l2}, similar to Case A and B.
If $u=4$, we may relabel the vertices on the quiver, $1\mapsto 1$ and $2+t\mapsto n-t$ 
for all $t\geq 0$. After relabeling, we are in the setting of $u=n-2$.
If $u=n-1$ or $u=3$, by relabeling with vertex $u$ as vertex $1$, we are in Case A, 
so the lemma follows from Lemma \ref{casea}.
\end{proof} 

\begin{example}
In this example we illustrate the construction in Case C. Consider the following quiver of Case C and $d=(1, \,1, \,2, \,2, \,1, \,2)$.
$$
\xymatrix@=5mm{
&& 1 \ar[d]\\
2 &\ar[l] 3 &\ar[l]4 \ar[r] &5\ar[r]^{\alpha}&6}
$$

(1) We first construct a rigid $A$-projective $D$-module $M^1\oplus M^2$ of $\Delta$-dimension vector $d'=(1, \,1, \,2, \,2, \,0, \,0)$
as follows.
$$
\xymatrix@=4mm{
M^1: &3\ar@{=>}[dr]\ar[dl] &&&&& M^2: & 2\ar@{=>}[dr]\\
2\ar@{=>}[dr] && 4=W_3\ar[dr]\ar[dl]&&&&&&3\ar@{=>}[dr]\ar[dl]\\
&3\ar[dl] && 5 \ar[dr] &&&& 2 \ar@{=>}[dr]&& 4=W_1\ar@{=>}[dr]\ar[dl] \ar[drr]\\
2 &&& &6 &&&& 3\ar[dl]\ar@{=>}[dr]&&1\ar[dl] &5\ar@{=>}[dll]\ar[dr]\\
&&& & &&& 2\ar@{=>}[dr] &&4 =W_2\ar[drr]\ar[dl]&& &6\ar@{=>}[dl]\\
&&& & &&&  &3\ar[dl]& &&5 \ar[dr]&\\
&&& & &&& 2 && &&&6
}
$$
Here $M^1<_4M^2$ and $\mathrm{Hom}(M^1, M^2)^0_4=0=\mathrm{Hom}(M^2, M^1)^0_4$. The
spaces $W_1, \, W_2,\,W_3 $ are given in the picture, $W=W_1\oplus W_3$ and
$H_W=\{ \left(\begin{matrix} a_1 & 0\\ 0&a_2\end{matrix}\right)\mid a_1, a_2 \in k^*\}$.

(2) A rigid $A$-projective $D$-module  $R\oplus T$  of $\Delta$-dimension vector $d''=( 0, \,0, \,0, \,0, \,1, \,2)$ is
$$\xymatrix@=2mm{
& 6=V_1\ar@{=>}[dl] \\
5\ar[dr]&&\oplus  & 6=V_3\\
& 6=V_2,
}$$
with $V_1, \, V_2, \, V_3$ are indicated in the picture. We have $V=V_2\oplus V_3$ and
$H_V=\{ \left(\begin{matrix} g_1 & x\\ 0&g_2\end{matrix}\right)\mid g_1, g_2\in k^*, x \in k\}$.

(3) Now the quiver $Q'$ is
$$
\xymatrix@=5mm{
Q': && 1 \ar[d]\\
2 \ar[r]&\ar[r] 3 &4  &5\ar[l]&6\ar[l]}
$$

and the representations  $Z(d')$ and $Z(d'')$ of  the quiver $Q'$ are
$$\xymatrix@=5mm{ 1 \ar[d]\\ 4 & \oplus & 3\ar[r] &4  & \mbox{and} & 5 &\ar[l] 6 &\oplus & 5.} $$
Clearly, $\mathrm{Aut}_{kQ'}Z(d')\cong H_W$ and $\mathrm{Aut}_{kQ'}Z(d'')\cong H_V$. A rigid representation of $Q'$ with
dimension vector $c(d)=\dimv Z(d')+\dimv Z(d'')=(1, \,0, \,1, \,2, \,2, \,1)$ is
$$
\xymatrix@=10mm{
&& k \ar[d]^{\tiny \left( \begin{matrix} 1\\0\end{matrix}\right)}\\
0\ar[r]&k \ar[r]_{\tiny \left( \begin{matrix} 0\\1\end{matrix}\right)} &k^2 &k^2\ar[l]^{\tiny z=\left( \begin{matrix} 1&0\\1&1\end{matrix}\right)}
&k\ar[l]^{ \tiny \left( \begin{matrix} 1\\0\end{matrix}\right)}
}
$$

(4) The matrix $z$ in (3) induces a rigid $A$-projective $D$-module $M$ of $\Delta$-dimension vector $d$ as follows.

$$
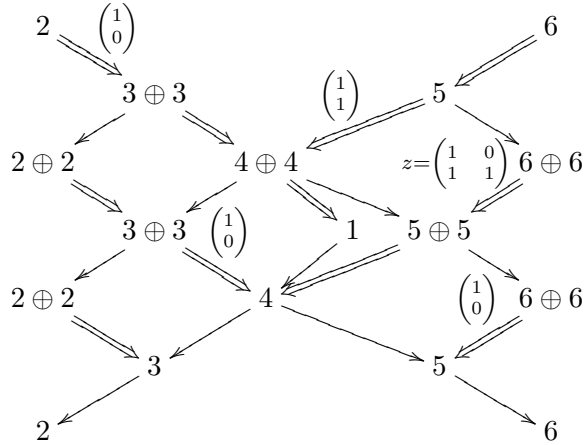
\xymatrix@=4mm{
  2\ar@{=>}[dr]^{\tiny \left( \begin{matrix} 1\\0\end{matrix}\right)}
&&&&& 6\ar@{=>}[dl]\\
&3\oplus 3\ar@{=>}[dr]\ar[dl] &&& 5\ar@{=>}[dll]_{\left( \tiny \begin{matrix} 1\\1\end{matrix}\right) }\ar[dr] \\
2\oplus 2 \ar@{=>}[dr]&& 4\oplus 4\ar@{=>}[dr]\ar[dl] \ar[drr] &&&6\oplus 6\ar@{=>}[dl]_{\tiny z=\left( \begin{matrix} 1&0\\1&1\end{matrix}\right)} \\
& 3\oplus 3\ar[dl]\ar@{=>}[dr]^{\left(\tiny \begin{matrix} 1\\0\end{matrix}\right)}&&1\ar[dl] &5\oplus 5\ar@{=>}[dll]\ar[dr]\\
2\oplus 2\ar@{=>}[dr] &&4\ar[drr]\ar[dl]&& &6\oplus 6\ar@{=>}[dl]_{\left(\tiny \begin{matrix} 1\\0\end{matrix}\right)}\\
&3\ar[dl]& &&5 \ar[dr]&\\
2 && &&&6
}
$$
As before, $i\oplus i$ means that the vector space is $2$-dimensional. If not specified, 
the matrices on the arrows connecting
two vector spaces of the same dimensions are identity matrices, and the 
other maps are determined by the commutativity relations.

We also compute directly the endomorphism ring to give another proof of the rigidity of $M$,
$$\End_DM=\{
f=(f_i)\mid (f_1, \, f_2, \, f_3, \, f_4, \, f_5|_{(M^1\oplus M^2)_5}, \, f_6|_{(M^1\oplus M^2)_6}) \in \End_D(M^1\oplus M^2),
$$$$(0, \,0, \, 0, \, 0, \, f_5|_{(R\oplus T)_5}, \, f_6|_{(R\oplus T)_6}) \in \End_D(R\oplus T) \mbox{ and } f_5M_{\beta^*}=M_{\beta^*}f_6
\}$$
and obtain $\dim_k \End_DM=15$, which is $\sum_{i=1}^6 d_i^2$, and so $M$ is rigid, by Theorem \ref{charlemma}.
\end{example}

\subsection{Case D} 

We have proven Theorem \ref{technicaltheorem} for Case A, B and C, and it only remains to prove for Case D. By relabeling
we are in Case A as follows.

\begin{lemma} \label{cased}
Let $Q$ be as in Case D and let $P$ be a projective representation of $Q$. 
Then there is a dimension vector $c$ such that the $\Aut_A P$-space $\Rep (D,P)$ 
is generically equivalent to the $\Gl(c)$-space $\Rep (Q',c)$ for some $c\in \mathbb{N}^n$.
\end{lemma}
\begin{proof}
We relabel the vertices on the quiver, $1\mapsto 1$ and $2+t\mapsto n-t$ 
for all $t\geq 0$. We are then in Case A, and so
the lemma follows from Lemma \ref{casea}.
\end{proof}

{\parindent=0cm
Bernt Tore Jensen, \\
Gj\o vik University College, \\
Postboks 191, \\
2802 Gj\o vik, \\
Norway, \\ 
email: bernt.jensen@hig.no \\\\}

{\parindent=0cm
Xiuping Su, \\
Department of Mathematical Sciences, \\
University of Bath,\\
Bath BA2 7JY,\\
United Kingdom.\\
email: xs214@bath.ac.uk }

\end{document}